\documentclass[12pt]{amsart}

\usepackage[all]{xy}

\newtheorem{theorem}{Theorem}[section]
\newtheorem{corollary}[theorem]{Corollary}
\newtheorem{lemma}[theorem]{Lemma}

\theoremstyle{remark}
\newtheorem{remark}{Remark}
\newtheorem*{claim}{Claim}
\newtheorem*{reason}{Proof}

\theoremstyle{definition}
\newtheorem{definition}[theorem]{Definition}
\newtheorem{example}[theorem]{Example}
\newtheorem{conditions}[theorem]{Conditions}

\newcounter{count}

\newcommand{\zerr}[2]{\ensuremath{\left(\!\!\!\begin{array}{c} #1\\ #2
\end{array}\!\!\!\right)}}

\begin{document}

\title{Partial Dynamical Systems and AF C$^*$-algebras}
\author{Justin R. Peters}
\address{Department of Mathematics, Iowa State University, Ames,
IA, USA} \email{peters@iastate.edu, ryzerr@iastate.edu}
\author{Ryan J. Zerr}

\begin{abstract}
We obtain a characterization in terms of dynamical systems of
those $r$-discrete groupoids for which the groupoid C$^*$-algebra
is approximately finite-dimensional (AF). These ideas are then
used to compute the $K$-theory for AF algebras by utilizing the
actions of these partial homeomorphisms, and these $K$-theoretic
calculations are applied to some specific examples of AF algebras.
Finally, we show that, for a certain class of dimension groups, a
groupoid can be obtained directly from the dimension group's
structure whose associated C$^*$-algebra has $K_0$ group
isomorphic to the original dimension group.
\end{abstract}

\bibliographystyle{hamsplain}

\maketitle

\section{Introduction}
\label{section1}

The use of dynamical systems to study C$^*$-algebras has received
much attention. In particular, ~\cite{tomiyama} discusses the
C$^*$-algebraic interpretations of a variety of dynamical systems
concepts and ~\cite{giordanoputnamskau, pimsnervoiculescu2,
pimsner} have considered, at least in part, AF subalgebras of
crossed product C$^*$-algebras and used the properties of the
dynamical system associated with the crossed product to gain
information about these algebras. In addition to this, the theory
of nonselfadjoint subalgebras of AF C$^*$-algebras has direct
connections with dynamical systems.

By a partial homeomorphism $\phi$ on a metric space $X$ we mean
that $\mbox{Dom}(\phi)$ and $\mbox{Ran}(\phi)$ are open subsets of
$X$ and $\phi$ is a homeomorphism from $\mbox{Dom}(\phi)$ to
$\mbox{Ran}(\phi)$.  Suppose $X$ is a $0$-dimensional space and
\[
\mathcal{R}_\phi=\{(x,\phi^n(x)):x\in\mbox{Dom}(\phi^n),n\in
\mathbb{Z}\}
\]
is the associated groupoid ($\phi^0=\mbox{id}|_X$).  Then
$\mathcal{R}_\phi$ is AF if and only if for any open set $U$
containing $\mbox{Dom}(\phi)$, $\bigcup_{n=0}^\infty \phi^n(U)
=X$~\cite[Corollary 4.7]{petersandpoon}.

The results of \S \ref{section2} arose from the question of how to
generalize this to the case where the groupoid is the union of the
graphs of a countable collection of partial homeomorphisms.  Note
that in the above result,
\[
\{\sigma_{V,n}=\phi^n|_V:n\in\mathbb{Z}, V\mbox{ a clopen set in
}\mbox{Dom}(\phi^n)\}
\]
is such a countable collection of partial homeomorphisms.  Thus,
the setting in \S \ref{section2} generalizes that
of~\cite{petersandpoon}.

In this paper we will consider partial homeomorphisms for which
the domain (and hence the range) are clopen subsets of the ambient
$0$-dimensional space. The main result of \S 2 is a dynamical
characterization of AF inverse semigroups (Definition
\ref{afsemigroup}).  We use this theorem (Theorem \ref{classify})
in Example \ref{odometer} to show that the inverse semigroup of
the odometer map is not AF.  Previously, the proof of this fact
had been via $K$-theory.

The importance of $K$-theory in the study of AF algebras makes it
natural to then consider, in \S \ref{section3}, how our dynamical
system characterization can be used to compute the $K_0$ group of
an AF algebra directly from a consideration of its groupoid.
Motivation for this is provided in~\cite{hermanandputnamandskau,
poon2, putnam}. In~\cite{putnam}
and~\cite{hermanandputnamandskau}, the $K$-theory of AF
subalgebras of certain crossed products is studied, and this
$K$-theory is computed through use of the short exact sequence of
Pimsner and Voiculescu~\cite{pimsnervoiculescu}. The computations
done here are more closely related to those of~\cite{poon}, but
valid for the more general partial dynamical systems considered as
part of the characterization in \S \ref{section2}.

The $K$-theoretic calculations of \S \ref{section3} are applied in
\S \ref{section4} to some specific examples of AF algebras,
including the CAR and GICAR algebras (Examples \ref{carexample}
and \ref{gicar}). It is here we find that, at least in some cases,
the $K_0$ group of an AF algebra can be realized as a group of
continuous functions acting on a $0$-dimensional compact Hausdorff
space.

For certain types of simple dimension groups, we discuss in \S
\ref{section5} a ``Pontryagin-like'' duality for groupoids which
uses the structure of the resulting character group to directly
construct an inverse semigroup of partial homeomorphisms and its
associated groupoid. The results of \S \ref{section3} can then be
used to prove that this groupoid is not only AF, but that the
$K_0$ group of the groupoid C$^*$-algebra is isomorphic to the
original dimension group.  This method brings full-circle our
ideas in \S \ref{section3} by reversing the process described
there. That is, whereas in \S \ref{section3} we compute the $K_0$
group directly from the groupoid, for the types of dimension
groups considered in \S \ref{section5}, we are able to compute the
groupoid directly from the $K_0$ group.

\section{AF Groupoids and Partial Dynamical Systems}
\label{section2}

In this section we will develop a characterization of those
groupoids that come from AF C$^*$-algebras. The conditions which
describe this characterization are dynamical in nature, and
therefore give a dynamical systems perspective on the nature of
these AF groupoids.

Given an AF algebra $\mathfrak{A}=\overline{\bigcup_{n\ge 1}
\mathfrak{A}_n}$, where $\{\mathfrak{A}_n\}_{n\ge 1}$ is an
increasing sequence of finite-dimensional subalgebras of
$\mathfrak{A}$, consider $\mathfrak{D}=\overline{\bigcup_{n\ge
1}\mathfrak{D}_n}$ where $\mathfrak{D}_n$ is the subalgebra of
$\mathfrak{A}_n$ spanned by the diagonal matrix units in
$\mathfrak{A}_n$.  As described in~\cite{power}, we can identify
the matrix units of each subalgebra $\mathfrak{A}_n$ with partial
homeomorphisms on the Gelfand space $X$ of $\mathfrak{D}$. The set
$X$ is a $0$-dimensional compact Hausdorff space, and the
equivalence relation $\mathcal{R}_{\mathfrak{A}}$ on $X$ that
results from taking the union of the graphs of these partial
homeomorphisms is an $r$-discrete groupoid and will be called the
groupoid of the AF algebra $\mathfrak{A}$.

In the theory of groupoid C$^*$-algebras,~\cite{paterson,
renault}, one forms the C$^*$-algebra of a groupoid $\mathcal{R}$,
$C^*(\mathcal{R})$, by considering the continuous functions of
compact support on $\mathcal{R}$.  In the case of the groupoid
$\mathcal{R}_{\mathfrak{A}}$, the C$^*$-algebra of the groupoid is
$\mathfrak{A}$ itself.  So, through a consideration of this
groupoid, we can establish necessary conditions on those groupoids
that give rise to AF algebras. The main result of this section
(Theorem \ref{classify}) will consider the converse. That is, we
will establish conditions that ensure a given groupoid gives rise
to a C$^*$-algebra which is AF.  This leads us to make the
following definition.

\begin{definition}
\label{afgroupoid}

A groupoid $\mathcal{R}$ will be called an {\it AF groupoid}
provided $C^*(\mathcal{R})$ is AF.

\end{definition}

To begin, we let $X$ be a $0$-dimensional compact Hausdorff space.
That is, $X$ is compact Hausdorff and has a basis consisting of
clopen sets.  Suppose that $\{\sigma_{r,s}^{(n)}: n\ge 0, 1\le
r\le m_n, 1\le s\le \kappa(r,n)\}$ is a countable collection of
partial homeomorphisms on $X$ where $\{m_n\}_{n=0}^\infty$ is a
sequence of positive integers (with $m_0=1$), and, for each $n\ge
1$ and $1\le r\le m_n$, $\kappa(r,n)\in\mathbb{Z}^+$.  Further
suppose that for any partial homeomorphism in this collection, say
$\sigma_{r,s}^{(n)}$, $\mbox{Dom}(\sigma_{r,s}^{(n)})$, the domain
of $\sigma_{r,s}^{(n)}$, is the clopen set $B(r,n)\subset X$. That
is to say, for $n$ and $r$ fixed, $\sigma_{r,1}^{(n)},\ldots,
\sigma_{r,\kappa(r,n)}^{(n)}$ share a common domain of $B(r,n)$.

We now form the clopen subsets $U_n$, for all $n\ge 1$, of $X$ by
\[
U_n=\bigcup_{r=1}^{m_n}B(r,n),
\]
and suppose that for each $n\ge 1$, the collection
$\{B(r,n)\}_{r=1}^{m_n}$ partitions $U_n$.  In particular, for
each fixed $n\ge 1$, the sets $\{B(r,n)\}_{r=1}^{m_n}$ are
pairwise disjoint.  For notational convenience, let
$U_0=B(1,0)=X$.

With the notation thus established, we assume that for each fixed
$n\ge 1$ the collection $\{\sigma_{r,s}^{(n)}: 1\le r\le m_n, 1\le
s\le\kappa(r,n)\}$ satisfies the following conditions.

\begin{conditions}
\label{conditions}
\end{conditions}

\begin{list}{(\roman{count})}{\usecounter{count}}
\item $\sigma_{r,1}^{(n)}(x)=x$, for all $1\le r\le m_{n}$ and
$x\in B(r,n)$; \item for each $1\le r\le m_n$ and $1\le s\le
\kappa(r,n)$, there exists an $i$, $1\le i\le m_{n-1}$, such that
$\sigma_{r,s}^{(n)}(B(r,n))\subset B(i,n-1)$; and \item
$\bigvee_{r=1}^{m_{n}}\bigvee_{s=1}^{\kappa(r,n)}
\sigma_{r,s}^{(n)}(B(r,n))$ is a partition of $U_{n-1}$.
\end{list}
We will consider the inverse semigroup generated by these partial
homeomorphisms (Definition \ref{isdefinition}).

\begin{remark}
From conditions (i) and (ii) it follows that the sequence
$\{U_n\}_{n=1}^\infty$ is nested decreasing.
\end{remark}

Now, given a collection $T$ of partial homeomorphisms on a
topological space $X$, let $A(T)$ be the set of partial
homeomorphisms on $X$ given by
\[
A(T)=\{\rho:\exists\sigma\in T \mbox{ and clopen subset }
B\subset\mbox{Dom}(\sigma)\mbox{ with }\rho=\sigma|_B\}.
\]
Note that $T\subset A(T)$.  We then make the following definition.

\begin{definition}
\label{isdefinition}

Given a collection $T$ of partial homeomorphisms on a topological
space $X$, the {\it inverse semigroup generated by $T$} will be
the smallest set of partial homeomorphisms on $X$ which contains
$A(T)$ and is closed for inverses and compositions (wherever those
compositions make sense).

\end{definition}

For the collection $\{\sigma_{r,s}^{(n)}\}$, we will now give a
description of this inverse semigroup. To begin, let $S$ be the
inverse semigroup generated by $\{\sigma_{r,s}^{(n)}\}$ and define
$\tau_{r,s}^{(1)}=\sigma_{r,s}^{(1)}$ for all $1\le r\le m_1$ and
$1\le s\le k(r,1)=\kappa(r,1)$.  By then taking the compositions
$\tau_{r_1,s_1}^{(1)}\circ\sigma_{r_2,s_2}^{(2)}$ of the partial
homeomorphisms $\{\sigma_{r,s}^{(2)}\}$ with the partial
homeomorphisms $\{\tau_{r,s}^{(1)}\}$, we obtain partial
homeomorphisms $\{\tau_{r',s'}^{(2)}\}$, for $1\le r'\le m_2$ and
$1\le s'\le k(r',s)$ for some $k(r',2)\in\mathbb{Z}^+$, which are
in $S$ and for which $\mbox{Dom}(\tau_{r',s'}^{(2)})=B(r',2)$ for
all $r'$ and $s'$. We proceed inductively to obtain, for all $n\ge
1$, the partial homeomorphisms $\{\tau_{r,s}^{(n)}\}$ where $1\le
r\le m_n$ and $1\le s\le k(r,n)$ for some $k(r,n)\in\mathbb{Z}^+$,
such that $\mbox{Dom}(\tau_{r,s}^{(n)}) =B(r,n)$, for all $r$ and
$s$, and where each $\tau_{r,s}^{(n)}$ is an element of $S$.

\begin{remark}
\label{remark2}
The collection $\{\tau_{r,s}^{(n)}: n\ge 1, 1\le
r\le m_n, 1\le s\le k(r,n)\}$ of partial homeomorphisms so
constructed is also a generating set for $S$ since it contains
$\sigma_{r,s}^{(n)}$ for every $n$, $r$, and $s$.
\end{remark}

\begin{remark}
\label{remark3}
For each $n\ge 1$,
\[
\bigvee_{r=1}^{m_n}\bigvee_{s=1}^{k(r,n)} \tau_{r,s}^{(n)}(B(r,n))
\]
forms a partition of $X$, and the sequence of partitions thus
obtained is a sequence of successively finer partitions of $X$.
\end{remark}

Now, let $\mathcal{R}$ be an $r$-discrete groupoid which is an
equivalence relation on $X$.  We will say that an inverse
semigroup of partial homeomorphisms $Z$ {\it generates}
$\mathcal{R}$ if the union of the graphs of the partial
homeomorphisms in $Z$ equals $\mathcal{R}$.  This leads us to make
the following definition.

\begin{definition}
\label{afsemigroup}

We will call an inverse semigroup $Z$ of partial homeomorphisms an
{\it AF inverse semigroup} provided the groupoid generated by $Z$
is AF.

\end{definition}

When a groupoid is generated by an inverse semigroup $Z$, we will
assume the groupoid has a topology with basis consisting of the
graphs of elements of $Z$ restricted to clopen subsets in $X$. For
the inverse semigroup $S$ constructed above and the groupoid
$\mathcal{R}$ generated by it, make note of the fact that the
collection
\[
\bigvee_{n=1}^\infty\bigvee_{r=1}^{m_n} \bigvee_{s,s'=1}^{k(r,n)}
\{\Gamma(\rho):\exists B\subset\mbox{Ran}(\tau_{r,s'}^{(n)}),
\mbox{clopen, with
}\rho=\tau_{r,s}^{(n)}\circ{\tau_{r,s'}^{(n)}}^{-1}|_B\}
\]
forms a basis for this topology on $\mathcal{R}$, where
$\Gamma(\rho)$ represents the graph of $\rho$ as a subset of
$X\times X$. We will show that any such groupoid is AF. Before
establishing this fact, however, we pause briefly to give a
heuristic account of why the converse is true.

The partial homeomorphisms
$\tau_{r,s}^{(n)}\circ{\tau_{r,s'}^{(n)}}^{-1}$ should be thought
of as playing a role corresponding to the matrix unit partial
homeomorphisms in the AF groupoid $\mathcal{R}_{\mathfrak{A}}$. In
particular, $\tau_{r,s}^{(n)} \circ{\tau_{r,s'}^{(n)}}^{-1}$ is
the partial homeomorphism corresponding to the matrix unit
$e_{s',s} \in M_{k(r,n)}$ where the AF algebra
$\mathfrak{A}=\overline{\bigcup_{n\ge 1} \mathfrak{A}_n}$ is such
that
\[
\mathfrak{A}_n=M_{k(1,n)}\oplus\cdots\oplus M_{k(m_n,n)}.
\]
Also, the sets $B(r,n)$ act in the same way as the supports of the
upper left diagonal matrix units in the subalgebras
$\mathfrak{A}_n$ when considering the groupoid of the AF algebra
$\mathfrak{A}$. With these motivating ideas in mind, it is fairly
straightforward to show that for each AF algebra $\mathfrak{A}$,
there exists a collection of partial homeomorphisms satisfying
Conditions \ref{conditions} which generate
$\mathcal{R}_{\mathfrak{A}}$ in the sense just described.
Therefore, Conditions \ref{conditions} describe necessary
conditions on AF groupoids.

In order to show that these conditions are sufficient we begin by
defining, for any $n\ge 1$, the set
\[
\mathcal{R}_n=\bigcup_{r=1}^{m_n}\bigcup_{s,s'=1}^{k(r,n)}
\Gamma\left(\tau_{r,s}^{(n)}\circ{\tau_{r,s'}^{(n)}}^{-1}\right)
\subset\mathcal{R},
\]
and give $\mathcal{R}_n$ the topology inherited from
$\mathcal{R}$.

\begin{lemma}
\label{groupoidstructure}

The sequence $\{\mathcal{R}_n\}_{n=1}^\infty$ is a nested
increasing sequence of groupoids such that
$\mathcal{R}=\bigcup_{n=1}^\infty \mathcal{R}_n$.

\end{lemma}

\begin{proof}

First, for any $n\ge 1$ let $(x,y)\in\mathcal{R}_n$.  That is to
say, there exists $r,s,s'$ such that $(x,y)\in\Gamma\left(
\tau_{r,s}^{(n)}\circ{\tau_{r,s'}^{(n)}}^{-1}\right)$.  So,
$y=\tau_{r,s}^{(n)}\circ{\tau_{r,s'}^{(n)}}^{-1}(x)$.  It follows
that
\[
{\tau_{r,s}^{(n)}}^{-1}(y)={\tau_{r,s'}^{(n)}}^{-1}(x)\in B(r,n).
\]
By construction, there exists $B(\overline{r},n+1)$ and
$\tau_{\overline{r},\overline{s}}^{(n+1)}$ such that
\[
{\tau_{r,s}^{(n)}}^{-1}(y)={\tau_{r,s'}^{(n)}}^{-1}(x)\in
\tau_{\overline{r},\overline{s}}^{(n+1)}(B(\overline{r},n+1)).
\]
So, for some $z\in B(\overline{r},n+1)$,
\[
\tau_{r,s}^{(n)}\circ\tau_{\overline{r},\overline{s}}^{(n+1)}(z)=y
\]
and
\[
\tau_{r,s'}^{(n)}\circ\tau_{\overline{r},\overline{s}}^{(n+1)}(z)=x.
\]
Again, by construction, there exist $s_1$ and $s_2$ such that
\[
\tau_{\overline{r},s_1}^{(n+1)}=\tau_{r,s}^{(n)}\circ\tau_{
\overline{r},\overline{s}}^{(n+1)}
\]
and
\[
\tau_{\overline{r},s_2}^{(n+1)}=\tau_{r,s'}^{(n)}\circ\tau_{
\overline{r},\overline{s}}^{(n+1)}.
\]
Hence, we see that $y=\tau_{\overline{r},s_1}^{(n+1)}(z)$ and
$z={\tau_{\overline{r},s_2}^{(n+1)}}^{-1}(x)$, which leaves
$y=\tau_{\overline{r},s_1}^{(n+1)}\circ{\tau_{\overline{r},s_2}^{
(n+1)}}^{-1}(x)$.  Therefore,
\[
(x,y)\in\bigcup_{r=1}^{m_{n+1}}\bigcup_{s,s'=1}^{k(r,n+1)}
\Gamma\left(\tau_{r,s}^{(n+1)}\circ{\tau_{r,s'}^{(n+1)}}^{-1}
\right)=\mathcal{R}_{n+1}.
\]
So, $\{\mathcal{R}_n\}_{n=1}^\infty$ is nested increasing.  Since
$\mathcal{R}=\bigcup_{n=1}^\infty \mathcal{R}_n$ is immediate, it
remains to show that each $\mathcal{R}_n$ is a groupoid.

To show $\mathcal{R}_n$ is a groupoid, the only nontrivial part is
to show that if $(x,y),(y,z)\in\mathcal{R}_n$ then
$(x,z)\in\mathcal{R}_n$, which we leave to the reader.

\end{proof}

\begin{remark}

By using arguments similar to those in Lemma
\ref{groupoidstructure}, one can show that if $B_1= \Gamma\left(
\tau_{r_1,s_1}^{(n_1)}\circ{\tau_{r_1,s_1'}^{(n_1)}}^{-1}\right)$
and $B_2=\Gamma\left( \tau_{r_2,s_2}^{(n_2)}\circ{
\tau_{r_2,s_2'}^{(n_2)}}^{-1}\right)$ are two basis elements of
$\mathcal{R}$ then either $B_1$ and $B_2$ are disjoint or one is a
subset of the other.  Hence, the topology of $\mathcal{R}$ has a
basis consisting of clopen sets.  It also follows easily that each
basis element is in fact compact.  Thus, since each
$\mathcal{R}_n$ is a finite union of compact open subsets of
$\mathcal{R}$, it follows that $\mathcal{R}_n$ is compact for all
$n\ge 1$.

\end{remark}

For $n\ge 1$ fixed, consider
$C_c(\mathcal{R}_n)=C(\mathcal{R}_n)$, the continuous functions
(of compact support) on $\mathcal{R}_n$.  We endow this set with
the groupoid convolution multiplication~\cite{paterson} in order
to obtain a $*$-algebra structure.  Let $1\le r\le m_n$, $g\in
C(B(r,n))$, and $1\le i,j\le k(r,n)$ be given, and define
$f_{i,j}^{(r)}\otimes g\in C(\mathcal{R}_n)$ by
\[
f_{i,j}^{(r)}\otimes g(x,y)=\left\{ \begin{array}{ll}
g\circ{\tau_{r,i}^{(n)}}^{-1}(x), & \mbox{if
$(x,y)\in\Gamma\left(\tau_{r,j}^{(n)}\circ{\tau_{
r,i}^{(n)}}^{-1}\right)$}\\
0, & \mbox{otherwise.}
\end{array}
\right.
\]
In the following lemma, which is somewhat similar to Lemma 2.2
of~\cite{poon}, we show that $C(\mathcal{R}_n)$ is spanned by
functions of this form.

\begin{lemma}
\label{isomorphic}

The algebra $C(\mathcal{R}_n)$ is isomorphic to
\[
\bigoplus_{r=1}^{m_n} M_{k(r,n)}\left( C(B(r,n))\right).
\]

\end{lemma}

\begin{proof}

First, we take $g\in C(\mathcal{R}_n)$.  Let
$\Gamma_{i,j}^{(r)}=\Gamma\left(
\tau_{r,j}^{(n)}\circ{\tau_{r,i}^{(n)}}^{-1}\right)$.  If we let
\[
g|_{\Gamma_{i,j}^{(r)}}(x,y)=\left\{ \begin{array}{ll} g(x,y), &
\mbox{if $(x,y)\in\Gamma_{i,j}^{(r)}$}\\
0, & \mbox{otherwise,}
\end{array}
\right.
\]
then it is clear that
\[
g=\sum_{r=1}^{m_n}\sum_{i,j=1}^{k(r,n)}g|_{\Gamma_{i,j}^{(r)}}.
\]
It is also straightforward to see that
$g|_{\Gamma_{i,j}^{(r)}}=f_{i,j}^{(r)}\otimes \tilde{g}$, with
\[
\tilde{g}=g\circ r^{-1}\circ\tau_{r,i}^{(n)} :B(r,n)\to\mathbb{C},
\]
where here $r:\Gamma_{i,j}^{(r)}\to\tau_{r,i}^{(n)} (B(r,n))$ is
projection onto the first coordinate and restricted to the subset
$\Gamma_{i,j}^{(r)}$ of $\mathcal{R}$. The map $r$ is essentially
the range map of Renault~\cite{renault}, and since both $r$ and
$\tau_{r,i}^{(n)}$ are homeomorphisms in this case, we see that
$\tilde{g}\in C(B(r,n))$. Hence, it follows that
$C(\mathcal{R}_n)$ is equal to the span of the collection
\[
\left\{ f_{i,j}^{(r)}\otimes g: 1\le r\le m_n, 1\le i,j\le k(r,n),
g\in C(B(r,n))\right\}.
\]

Now, for $1\le r_1<r_2\le m_n$, $g_1\in C(B(r_1,n))$, $g_2\in
C(B(r_2,n))$, $1\le i_1,j_1\le k(r_1,n)$, and $1\le i_2,j_2\le
k(r_2,n)$ we have
\[
\left( f_{i_1,j_1}^{(r_1)}\otimes g_1\right)\ast\left(
f_{i_2,j_2}^{(r_2)}\otimes g_2\right)=0
\]
since $\tau_{r_1,j_1}^{(n)}(B(r_1,n))
\cap\tau_{r_2,i_2}^{(n)}(B(r_2,n))=\emptyset$, where $\ast$ is the
convolution product on $C(\mathcal{R}_n)$.

Next, let $1\le r\le m_n$ be fixed and take $1\le i,j,l,m\le
k(r,n)$. If $j\ne l$ then
$\tau_{r,j}^{(n)}(B(r,n))\cap\tau_{r,l}^{(n)}(B(r,n))=\emptyset$,
and so, for $g_1,g_2\in C(B(r,n))$, we again have
\[
\left(f_{i,j}^{(r)}\otimes g_1\right)\ast\left(
f_{l,m}^{(r)}\otimes g_2\right) =0.
\]
However, if $j=l$ take $(x,z)\in\Gamma_{i,j}^{(r)}$ and
$(z,y)\in\Gamma_{l,m}^{(r)}$. Then,
\begin{eqnarray*}
\left(f_{i,j}^{(r)}\otimes g_1\right)\ast\left( f_{l,m}^{(r)}
\otimes g_2\right) (x,y)&=& f_{i,j}^{(r)} \otimes g_1(x,z) \cdot
f_{l,m}^{(r)}\otimes g_2(z,y)\\
&=& g_1\circ{\tau_{r,i}^{(n)}}^{-1}(x) \cdot g_2\circ
{\tau_{r,l}^{(n)}}^{-1}(z).
\end{eqnarray*}
Also, we see that $f_{i,m}^{(r)}\otimes(g_1\cdot
g_2)(x,y)=(g_1\cdot g_2)\circ {\tau_{r,i}^{(n)}}^{-1}(x)$ since
\begin{eqnarray*}
y&=&\tau_{r,m}^{(n)}\circ{\tau_{r,l}^{(n)}}^{-1}(z)
=\tau_{r,m}^{(n)}\circ{\tau_{r,l}^{(n)}}^{-1}\circ
\tau_{r,j}^{(n)}\circ {\tau_{r,i}^{(n)}}^{-1}(x)\\
&=&\tau_{r,m}^{(n)}\circ{\tau_{r,i}^{(n)}}^{-1}(x).
\end{eqnarray*}
So,
\[
f_{i,m}^{(r)}\otimes(g_1\cdot g_2)(x,y)=
g_1\circ{\tau_{r,i}^{(n)}}^{-1}(x)\cdot
g_2\circ{\tau_{r,i}^{(n)}}^{-1} (x).
\]
Since $z=\tau_{r,j}^{(n)}\circ{\tau_{r,i}^{(n)}}^{-1}(x)$ we see
that ${\tau_{r,i}^{(n)}}^{-1}(x)={\tau_{r,l}^{(n)}}^{-1} (z)$ due
to the fact that $j=l$.  Hence, for all $(x,y)\in\mathcal{R}_n$,
\[
f_{i,m}^{(r)}\otimes(g_1\cdot g_2)(x,y)=\left(
f_{i,j}^{(r)}\otimes g_1\right)\ast\left( f_{l,m}^{(r)}\otimes
g_2\right)(x,y).
\]
In particular,
\[
\left(f_{i,j}^{(r)}\otimes\chi_{_{B(r,n)}}\right)\ast\left(
f_{l,m}^{(r)}\otimes \chi_{_{B(r,n)}}\right)=\delta_{j,l}\left(
f_{i,m}^{(r)}\otimes\chi_{_{B(r,n)}}\right),
\]
where $\delta_{j,l}=1$ if $j=l$ and $0$ if $j\ne l$.

Next, for $r,g,i,j$ as above, we see that
\[
f_{j,i}^{(r)}\otimes\overline{g}(x,y)=\left\{
\begin{array}{ll}
\overline{g}\circ{\tau_{r,j}^{(n)}}^{-1}(x),&\mbox{if $(x,y)\in
\Gamma_{j,i}^{(r)}$}\\
0, &\mbox{otherwise}
\end{array}
\right.
\]
and
\[
f_{i,j}^{(r)}\otimes\overline{g}(y,x)=\left\{
\begin{array}{ll}
\overline{g}\circ{\tau_{r,i}^{(n)}}^{-1}(y),&\mbox{if $(y,x)\in
\Gamma_{i,j}^{(r)}$}\\
0, &\mbox{otherwise.}
\end{array}
\right.
\]
Thus, since $(x,y)\in\Gamma_{j,i}^{(r)}$ implies not only that
$(y,x)\in\Gamma_{i,j}^{(r)}$ but also that
${\tau_{r,j}^{(n)}}^{-1}(x)={\tau_{r,i}^{(n)}}^{-1}(y)$, we see
that in fact,
\[
f_{j,i}^{(r)}\otimes\overline{g}(x,y)= f_{i,j}^{(r)}\otimes
\overline{g}(y,x),
\]
for all $(x,y)\in\mathcal{R}_n$.  Therefore,
\begin{eqnarray*}
\left(f_{i,j}^{(r)}\otimes g\right)^*(x,y)&=& \overline{
f_{i,j}^{(r)} \otimes g}(y,x)
=f_{i,j}^{(r)}\otimes \overline{g}(y,x)\\
&=&f_{j,i}^{(r)}\otimes\overline{g}(x,y).
\end{eqnarray*}
So, $\left(f_{i,j}^{(r)}\otimes\chi_{_{B(r,n)}}\right)^*
=f_{j,i}^{(r)}\otimes \chi_{_{B(r,n)}}$.

Finally, because the images of the partial homeomorphisms
$\{\tau_{r,i}^{(n)}\}$ partition $X$, we have that
\[
\sum_{r=1}^{m_n}\sum_{i=1}^{k(r,n)}f_{i,i}^{(r)}\otimes
\chi_{_{B(r,n)}}= \chi_{_{\Delta(X_{\mathfrak{A}})}}
\]
is the identity on $C(\mathcal{R}_n)$.

This shows that for each $1\le r\le m_n$, the collection of
functions
$\{f_{i,j}^{(r)}\otimes\chi_{_{B(r,n)}}\}_{i,j=1}^{k(r,n)}\subset
C(\mathcal{R}_n)$ forms a system of matrix units for a full
$k(r,n)\times k(r,n)$ matrix algebra.  This completes the proof.

\end{proof}

Now, considering the sequence of algebras $\{C(\mathcal{R}_n)
\}_{n=1}^\infty$, define, for each $n\ge 1$, the map
$\phi_{n,n+1}:C(\mathcal{R}_n)\to C(\mathcal{R}_{n+1})$ by
\[
\phi_{n,n+1}(f)(x,y)=\left\{ \begin{array}{ll} f(x,y), & \mbox{if
$(x,y)\in\mathcal{R}_n$}\\
0, & \mbox{otherwise.}
\end{array}
\right.
\]
Since $\mathcal{R}_n$ is a subset of $\mathcal{R}_{n+1}$, it
follows that $\phi_{n,n+1}(f)$ is in fact an element of
$C(\mathcal{R}_{n+1})$. It is also fairly clear that
$\phi_{n,n+1}$ is a unital injective $*$-homomorphism. We can
therefore view $C(\mathcal{R}_n)$ as a subalgebra of
$C(\mathcal{R}_{n+1})$, and in fact, via the map $\Phi_n:
C(\mathcal{R}_n)\to C_c(\mathcal{R})$, where
\[
\Phi_{n}(f)(x,y)=\left\{ \begin{array}{ll} f(x,y), & \mbox{if
$(x,y)\in\mathcal{R}_n$}\\
0, & \mbox{otherwise,}
\end{array}
\right.
\]
we embed $C(\mathcal{R}_n)$ as a $*$-subalgebra of
$C_c(\mathcal{R})$.  As such, we have
\[
\bigcup_{n=1}^\infty C(\mathcal{R}_n)\subset C_c(\mathcal{R}).
\]
If, however, we let $f\in C_c(\mathcal{R})$, then by the
compactness of $\mbox{supp}(f)$, the support of $f$, there exists
an $N\ge 1$ such that
\[
\mbox{supp}(f)\subset\bigcup_{r=1}^{m_N}\bigcup_{s,s'=1}^{k(r,N)}
\Gamma\left(\tau_{r,s}^{(n)}\circ{\tau_{r,s'}^{(n)}}^{-1}\right).
\]
So, $f|_{\mathcal{R}_N}\in C(\mathcal{R}_N)$ and we see that
$\Phi_N\left(f|_{\mathcal{R}_N}\right)=f$.  Therefore,
\[
\bigcup_{n=1}^\infty C(\mathcal{R}_n)=C_c (\mathcal{R}).
\]

Now, if $W$ is a $0$-dimensional compact Hausdorff space, there
exists a sequence $\{E_n\}_{n=1}^\infty$ of successively finer
partitions of $W$ which generate the topology.  Thus, $C(W)$ is
the closed union of an increasing sequence of subalgebras,
$\{C(E_n)\}_{n=1}^\infty$, where $C(E_n)$ consists of those
functions that are constant on elements of the partition $E_n$.
Hence, $C(W)$ is AF. It follows that for $k\ge 1$, $M_k(C(W))$ is
AF. Therefore, for all $n\ge 1$, $C(\mathcal{R}_n)$ is an AF
algebra, and so, by~\cite{bratteli},
\[
C^*(\mathcal{R})=\overline{C_c(\mathcal{R})}
\]
is an AF algebra.  We summarize in the following theorem.

\begin{theorem}
\label{classify}

Suppose there exists a nested decreasing sequence of clopen
subsets $\{U_n\}_{n=0}^\infty$ of $X$
\textup{(}$0$-dimensional\textup{)}, clopen sets
$\{B(r,n)\}_{r=1}^{m_n}$ for each $n\ge 0$ which partition $U_n$
\textup{(}with $U_0=B(1,0)=X$\textup{)}, and a countable
collection
\[
\{\sigma_{r,s}^{(n)}:n\ge 0,1\le r\le m_n,1\le s\le \kappa(r,n)\}
\]
of partial homeomorphisms on $X$ such that $\mbox{Dom}(
\sigma_{r,s}^{(n)})=B(r,n)$ for all $n$, $r$, and $s$, and which
satisfy Conditions \textup{\ref{conditions}}.  Then, the
$r$-discrete groupoid corresponding to the inverse semigroup
generated by this collection of partial homeomorphisms is AF.

\end{theorem}

This theorem, in conjunction with our earlier remarks leads to the
following corollaries, the first of which is a generalization
of~\cite[Corollary $4.7$]{petersandpoon}.

\begin{corollary}
\label{classifycor}

The $r$-discrete groupoid $\mathcal{R}$ is AF if and only if it is
an equivalence relation on the $0$-dimensional space $X$ and there
exists a countable collection of partial homeomorphisms
\textup{(}as in Theorem \textup{\ref{classify})} satisfying
Conditions \textup{\ref{conditions}} which generate $\mathcal{R}$.

\end{corollary}

\begin{corollary}
\label{classifycor2}

An inverse semigroup of partial homeomorphisms on $X$
\textup{(}$0$-dimensional\textup{)} is AF if and only if it is
generated by a subset satisfying Conditions
\textup{\ref{conditions}}.

\end{corollary}

\begin{remark}

We here briefly mention a connection between our results and the
theory of nonselfadjoint subalgebras of C$^*$-algebras.
Specifically, suppose $T$ is a set of partial homeomorphisms on
$X$ and consider the inverse semigroup $Z$ generated by $T$.
Suppose $V \subset Z$ is a set of partial homeomorphisms
containing $A(T)$ and for which $V\circ V\subset V$,$V \cup
{V}^{-1}=Z$, and
\[
V \cap{V}^{-1}=\left\{ \rho:\mbox{ there exists a clopen subset
}B\subset X\mbox{ with }\rho=\mbox{id}|_B\right\}.
\]
In the case that $T$ is a countable collection of partial
homeomorphisms satisfying Conditions \ref{conditions}, Corollary
\ref{classifycor2} tells us that $Z$ is AF. If $\mathcal{R}'$ is
the groupoid generated by $Z$, then the set $V$ will be such that
its associated nonselfadjoint subalgebra of $C^*(\mathcal{R}')$ is
a triangular AF algebra~\cite{peterspoonandwagner}. In the context
of the above discussion, $T$ can be taken to be the collection
\[
\{\sigma_{r,s}^{(n)}:n\ge 0,1\le r\le m_n,1\le s\le \kappa(r,n)\}
\]
of partial homeomorphisms and the inverse semigroup $Z$ will then
correspond to $S$.

\end{remark}

We will also state here a lemma that, in addition to being useful
in the following example, is of interest in its own right.

\begin{lemma}
\label{examplecorollary}

Suppose $Z$ is an AF inverse semigroup of partial homeomorphisms
on $X$.  If $\rho\in Z$ and $\rho\left(
\mbox{Dom}(\rho)\right)=\mbox{Dom}(\rho)$ then $\rho=\mbox{id}|_{
\mbox{Dom}(\rho)}$.

\end{lemma}

\begin{proof}

By Corollary \ref{classifycor2}, there exists a generating set for
$Z$ which satisfies Conditions \ref{conditions}. Thus, Remark
\ref{remark2} implies that if $\rho\in Z$ then there exists an
$n\ge 0$, $1\le r\le m_n$, and $1\le s,s'\le k(r,n)$ such that
\[
\rho=\tau_{r,s}^{(n)}\circ{\tau_{r,s'}^{(n)}}^{-1}|_B
\]
for some clopen subset $B\subset X$.  Now, by Remark
\ref{remark3}, if $s\ne s'$ then $\tau_{r,s}^{(n)}(B(r,n))
\cap\tau_{r,s'}^{(n)}(B(r,n))=\emptyset$.  Thus, $\rho(B)\cap
B=\emptyset$.  Hence, $\rho(B)=B$ implies $s=s'$, and so,
$\rho=\mbox{id}|_B$.

\end{proof}

\begin{example}
\label{odometer}

To illustrate the use of Theorem \ref{classify}, we use it to
prove a groupoid is not AF.  In particular, let
$X=\prod_{n=1}^\infty \{0,1\}_n$ be the set of all sequences of
$0$'s and $1$'s with the product topology.  Define the odometer
map $\phi:X\to X$ by $\phi((x_n))=(y_n)$, where
\[
y_n=\left\{ \begin{array}{lll} 0, &
\mbox{if $x_i=1$ for all $1\le i\le n$}\\
1, & \mbox{if $x_i=1$ for all $1\le i\le n-1$ and $x_n=0$}\\
x_n, & \mbox{if $x_i=0$ for some $1\le i\le n-1$.}
\end{array}
\right.
\]
This map is a homeomorphism on $X$ and we will show that the
groupoid $\overline{\mathcal{R}}=\{(x,\phi^n(x)):x\in
X,n\in\mathbb{Z}\}$ is not AF.

Considering $\phi$, we let $B=\{(x_n)\in X:x_1=a_1\}$, for some
fixed $a_1$. Then, $\phi^{2}(B)=B$ but, since
\[
\phi^{2}((a_1,0,0,\ldots))=(a_1,1,0,0,\ldots),
\]
we see that $\phi^{2}|_B\ne \mbox{id}|_B$.  Hence, the inverse
semigroup which generates $\overline{\mathcal{R}}$ is not AF by
Lemma \ref{examplecorollary}.

\begin{remark}

The groupoid $\overline{\mathcal{R}}$ is such that
$C^*(\overline{\mathcal{R}})$ is the crossed product C$^*$-algebra
$C(X)\times_\phi\mathbb{Z}$. Hence, the result of this example
could be obtained by appealing to $K$-theory and showing that
$K_1(C^*(\overline{\mathcal{R}}))\ne 0$. Since,
cf.~\cite{davidson}, the $K_1$ group of an AF algebra must be $0$,
it follows that $C^*(\overline{\mathcal{R}})$ is not AF.

\end{remark}

\end{example}

\section{AF Groupoids and $K$-theory}
\label{section3}

As we know from Elliott's Theorem~\cite{elliott}, the dimension
group of an AF algebra is a complete isomorphism invariant for the
algebra. In this section we use the partial dynamical systems
described in the previous section to calculate the dimension group
of the algebra.

To begin, we recall some facts from $K$-theory, of
which~\cite{davidson} and~\cite{blackadar} are good references. In
a C$^*$-algebra $\mathfrak{A}$, the projections $P$ and $Q$ are
said to be $\ast$-equivalent if there exists an element $R$ in the
algebra such that $P=R^*R$ and $Q=RR^*$.  Two projections in
$\bigcup_{n\ge 1}M_n(\mathfrak{A})$, say $P\in M_m(\mathfrak{A})$
and $Q\in M_n(\mathfrak{A})$ with $m\le n$, are then said to be
equivalent if $P\oplus 0_{n-m}$ is $\ast$-equivalent to $Q$ in
$M_n(\mathfrak{A})$. This is an equivalence relation on the set of
all projections in $\bigcup_{n\ge 1}M_n(\mathfrak{A})$ when
$\mathfrak{A}$ is an AF algebra.  We let $K_0^+(\mathfrak{A})$
denote the collection of all such equivalence classes.  The $K_0$
group of $\mathfrak{A}$, $K_0(\mathfrak{A})$, is then defined to
be the Grothendieck group of $K_0^+(\mathfrak{A})$.  That is,
$K_0(\mathfrak{A})$ is the set of all formal differences of
elements in $K_0^+(\mathfrak{A})$ modulo the identification of
$a_1-b_1$ with $a_2-b_2$ if and only if $a_1+b_2=a_2+b_1$.

To describe the dimension group of the algebra $\mathfrak{A}$, one
also needs the order unit of $K_0(\mathfrak{A})$.  In particular,
if $\mathfrak{A}$ is unital, the order unit $u$ is simply $u=[I]$,
the equivalence class of the identity.  Then, the dimension group
of the AF algebra $\mathfrak{A}$ is the scaled ordered group
\[
(K_0(\mathfrak{A}),K_0^+(\mathfrak{A}),u).
\]

As mentioned in the previous section, if $Z$ is compact Hausdorff
with a countable basis consisting of clopen sets, then there
exists a sequence $\{E_n\}_{n=1}^\infty$ of successively finer
partitions of $Z$ which generate the topology on $Z$.  Thus,
$C(Z)$ is the closed union of an increasing sequence of
subalgebras, $\{C(E_n)\}_{n=1}^\infty$, where $C(E_n)$ consists of
those functions that are constant on elements of the partition
$E_n$.  Hence, $C(Z)$ is AF.  It follows from this
and~\cite{davidson} that $K_0(C(Z))$ is the direct limit of the
sequence $\{K_0(C(E_n))\}=\{C(E_n,\mathbb{Z})\}$, and so,
\[
K_0(C(Z))=C(Z,\mathbb{Z}),
\]
with $K_0^+(C(Z))=\{g\in C(Z,\mathbb{Z}):g\ge 0\}$ (see also
Example \ref{continuousoncantor}).

Now, suppose we have the AF groupoid $\mathcal{R}$, which, as in
the notation of the previous section, is the union of the
increasing sequence $\{\mathcal{R}_n\}_{n=1}^\infty$.  We will use
the following lemma to calculate the dimension group of the AF
algebra $\mathfrak{A}$ corresponding to $\mathcal{R}$.

\begin{lemma}

For each $n\ge 1$,
\[
K_0(C(\mathcal{R}_n))\cong C(U_n,\mathbb{Z})
\]
and $K_0^+(C(\mathcal{R}_n))=\{g\in C(U_n,\mathbb{Z}):g\ge 0\}$.

\end{lemma}

\begin{proof}

We have by Lemma \ref{isomorphic} that
\begin{eqnarray*}
K_0(C(\mathcal{R}_n))&\cong & K_0\left(
\bigoplus_{r=1}^{m_n} M_{k(r,n)}(C(B(r,n)))\right)\\
&\cong&
\bigoplus_{r=1}^{m_n}K_0\left(M_{k(r,n)}(C(B(r,n)))\right).
\end{eqnarray*}
Since $K_0\left(M_{k(r,n)}(C(B(r,n)))\right)=K_0(C(B(r,n)))$, we
then have
\begin{eqnarray*}
K_0(C(\mathcal{R}_n))&\cong&\bigoplus_{r=1}^{m_n}K_0(C(B(r,n)))\\
&\cong& \bigoplus_{r=1}^{m_n}C(B(r,n),\mathbb{Z}).
\end{eqnarray*}
Because $U_n$ is partitioned by $\{B(r,n)\}_{r=1}^{m_n}$, we have
the desired result.

\end{proof}

Since for $Z$ $0$-dimensional, $C(Z)$ is AF, it follows that for
$k\ge 1$, $M_k(C(Z))$ is AF, and in fact, using our earlier
notation, that
\[
M_k(C(Z))=\overline{\bigcup_{n\ge 1}M_k(C(E_n))}.
\]
Furthermore, we note that $M_k(C(E_n))\cong
\oplus_{l=1}^{|E_n|}M_k$. From Davidson~\cite{davidson}, any
projection $P\in M_k(C(Z))$ is $\ast$-equivalent to a projection
$Q$ in $M_k(C(E_n))$ for some $n$.  But, a standard result from
matrix theory then tells us that $Q$ is unitarily equivalent to a
diagonal matrix in $\oplus_{l=1}^{|E_n|}M_k$.  It is an easy
result, cf.~\cite{blackadar}, which will not be repeated here,
that if $A$ and $B$ are unitarily equivalent projections in a
unital C$^*$-algebra then $A$ and $B$ are $\ast$-equivalent. So,
by this result, $P$ is $\ast$-equivalent to a diagonal projection
in $\oplus_{l=1}^{|E_n|}M_k$.  Hence, every projection in
$M_k(C(Z))$ is equivalent to a diagonal element
$\mbox{diag}(\lambda_1,\ldots,\lambda_k)$ of $M_k(C(Z))$, where
each of the $\lambda_1,\ldots,\lambda_k$ is a characteristic
function of a clopen subset of $Z$.

Now, taking a projection $p\in
C(\mathcal{R}_n)\cong\oplus_{r=1}^{m_n}M_{k(r,n)}(C(B(r,n)))$, we
will determine $[p]$, the equivalence class of $p$.  By the
arguments immediately above, we may assume without a loss of
generality that $p$ is a direct sum of diagonal projections in
$\oplus_{r=1}^{m_n} M_{k(r,n)}(C(B(r,n)))$.  So,
\[
p=\sum_{r=1}^{m_n}\sum_{i=1}^{k(r,n)}f_{i,i}^{(r)}\otimes
\chi_{_{A_{r,i}}}
\]
where the $A_{r,i}$ are clopen subsets of $B(r,n)$, for all $r$
and $i$.  Thus,
\[
[p]=\sum_{r=1}^{m_n}\sum_{i=1}^{k(r,n)}\chi_{_{A_{r,i}}}.
\]
In particular, this gives the order unit $u_n$ of
$K_0(C(\mathcal{R}_n))$ as
\begin{eqnarray*}
u_n&=&\left[\chi_{_{\Delta(X)}}\right]=\left[\sum_{r=1}^{m_n}
\sum_{i=1}^{k(r,n)} f_{i,i}^{(r)}\otimes\chi_{_{B(r,n)}}\right]\\
&=&\sum_{r=1}^{m_n}\sum_{i=1}^{k(r,n)}\chi_{_{B(r,n)}}\\
&=&\sum_{r=1}^{m_n} k(r,n)\chi_{_{B(r,n)}}.
\end{eqnarray*}
We summarize in the following lemma.

\begin{lemma}

For each $n\ge 1$ we have $K_0(C(\mathcal{R}_n))\cong
C(U_n,\mathbb{Z})$, with
\[
K_0^+(C(\mathcal{R}_n))=\{g\in C(U_n,\mathbb{Z}):g\ge 0\}
\]
and order unit $u_n=\sum_{r=1}^{m_n}k(r,n)\chi_{_{B(r,n)}}$.

\end{lemma}

This lemma, along with the fact that $C^*(\mathcal{R})$ is the
closure of the union of the increasing sequence
$\{C(\mathcal{R}_n)\}$, gives us all we need to calculate the
dimension group of $C^*(\mathcal{R})$.  We do so in the following
theorem, which is similar to~\cite[Corollary $3.3$]{poon}.

\begin{theorem}
\label{dimgroup}

The dimension group of $C^*(\mathcal{R})$ is equal to the direct
limit of the scaled ordered groups $K_0(C(\mathcal{R}_n))\cong
C(U_n,\mathbb{Z})$ where the order unit $u$ is the direct limit of
the order units
\[
u_n=\sum_{r=1}^{m_n}k(r,n)\chi_{_{B(r,n)}}
\]
and the connecting homomorphisms
$\phi_{(n,n+1)\ast}:C(U_n,\mathbb{Z})\to C(U_{n+1},\mathbb{Z})$
are given by
\[
\phi_{(n,n+1)\ast}(\xi)=\sum_{r=1}^{m_n}\sum_{l=1}^{t(r)}
\widehat{\xi\circ\tau_{r_l,s_l}^{(n+1)}}
\]
where
\[
B(r,n)=\bigcup_{l=1}^{t(r)}\tau_{r_l,s_l}^{(n+1)} (B(r_l,n+1))
\]
and $\widehat{\xi\circ\tau_{r_l,s_l}^{(n+1)}}$ is the function
\[
\widehat{\xi\circ\tau_{r_l,s_l}^{(n+1)}}(x)=\left\{
\begin{array}{ll} \xi\circ\tau_{r_l,s_l}^{(n+1)}(x), &
\mbox{if $x\in B(r_l,n+1)$}\\
0, & \mbox{otherwise.}
\end{array}
\right.
\]

\end{theorem}

\begin{proof}

We have proven everything except for the nature of the connecting
homomorphisms.  Let $n\ge 1$ be given and take $\xi\in
C(U_n,\mathbb{Z})$, $0\le\xi\le u_n$. Define, for $1\le r\le m_n$
and $1\le i\le k(r,n)$, the sets
\[
A_{r,i}=B(r,n)\cap\xi^{-1}([i,k(r,n)]),
\]
where $[i,k(r,n)]$ is the set of integers $\{i,i+1,\ldots,
k(r,n)\}$. These sets are clopen by the continuity of $\xi$. Let
\[
p=\sum_{r=1}^{m_n}\sum_{i=1}^{k(r,n)}f_{i,i}^{(r)}\otimes\chi_{_{
A_{r,i}}}\in C(\mathcal{R}_n).
\]
This is clearly a projection for which $[p]=\xi$.  It follows from
considering $\phi_{(n,n+1)}(p)$ that
\[
\phi_{(n,n+1)\ast}(\xi)=\sum_{r=1}^{m_n}\sum_{l=1}^{t(r)}\widehat{\xi\circ
\tau_{r_l, s_l}^{(n+1)}}.
\]

To show that this formula holds for an arbitrary element of
$C(U_n,\mathbb{Z})$, we use the fact that $\{\xi\in
C(U_n,\mathbb{Z}):0\le\xi\le u_n\}$ is a scale for
$C(U_n,\mathbb{Z})$.  Thus, every element of $C(U_n,\mathbb{Z})$
can be written as a linear combination $\sum_j\alpha_j\xi_j$ with
$0\le\xi_j\le u_n$, for all $j$. This completes the proof.

\end{proof}

\section{Examples}
\label{section4}

We will now investigate some examples where Theorem \ref{dimgroup}
is used to calculate the dimension group of specific AF groupoids.
To simplify the calculations, however, we will use the following
lemma.

\begin{lemma}
\label{dimgrplemma}

Retaining the notation from Theorem $\ref{dimgroup}$, let
\[
\widetilde{C}(U_n,\mathbb{Z})=\{f\in C(U_n,\mathbb{Z}):
f|_{B(r,n)}=\alpha(r,n)\chi_{_{B(r,n)}}, 1\le r\le m_n\}.
\]
If the collection of sets
\[
\bigvee_{n=1}^\infty\bigvee_{r=1}^{m_n}\bigvee_{s=1}^{k(r,n)}
\tau_{r,s}^{(n)}(B(r,n))
\]
is a basis for the topology on $X$, then
\[
\lim_{\longrightarrow} C(U_n,\mathbb{Z})\cong
\lim_{\longrightarrow} \widetilde{C}(U_n,\mathbb{Z}),
\]
where the connecting maps $\widetilde{\phi}_{(n,n+1)\ast}:
\widetilde{C}(U_n,\mathbb{Z})\to
\widetilde{C}(U_{n+1},\mathbb{Z})$ are just the restrictions of
the maps $\phi_{(n,n+1)\ast}$ to the sets
$\widetilde{C}(U_n,\mathbb{Z})$.

\end{lemma}

\begin{proof}

To begin, let $n\ge 1$ be fixed and choose $f\in
C(U_n,\mathbb{Z})$.  By the continuity of $f$ and the compactness
of $U_n$, there exists a partition $\{\mathcal{O}_1,\ldots,
\mathcal{O}_k\}$ of $U_n$ such that $f=\sum_{i=1}^k
\alpha_i\chi_{_{\mathcal{O}_i}}$.  Since we are assuming that
\[
\bigvee_{n=1}^\infty\bigvee_{r=1}^{m_n}\bigvee_{s=1}^{k(r,n)}
\tau_{r,s}^{(n)}(B(r,n))
\]
is a basis for the topology on $X$, we may assume, by refining the
partition $\{\mathcal{O}_1,\ldots, \mathcal{O}_k\}$ if necessary,
that for each $1\le i\le k$,
\[
\mathcal{O}_i=\tau_{r_{l(1)},s_{l(1)}}^{(n+1)}\circ\cdots\circ
\tau_{r_{l(m)},s_{l(m)}}^{(n+m)}(B(r_{l(m)},n+m)),
\]
where $m$ is a constant not depending on $i$.  Thus, by Theorem
\ref{dimgroup}, $\phi_{(n,n+1)\ast}(f)=\sum_{i=1}^{\overline{k}}
\beta_i\chi_{_{\overline{\mathcal{O}}_i}}$, where each set
$\overline{\mathcal{O}}_i$ is of the form
\[
\overline{\mathcal{O}}_i=\tau_{r_{l(2)},s_{l(2)}}^{(n+2)}\circ
\cdots\circ\tau_{r_{l(m)},s_{l(m)}}^{(n+m)}(B(r_{l(m)},n+m)).
\]
Therefore, for each $n\ge 1$ and $f\in C(U_n,\mathbb{Z})$, there
exists an $N\ge n$ such that
\[
\phi_{(N-1,N)\ast}\circ\cdots\circ\phi_{(n,n+1)\ast}(f)\in
\widetilde{C} (U_N,\mathbb{Z}).
\]
The result now follows from this by applying the universal
property of direct limits.

\end{proof}

\begin{example}
\label{carexample}

Consider the CAR algebra $\mathfrak{A}=\overline{\bigcup_{n\ge 1}
\mathfrak{A}_n}$, where $\mathfrak{A}_n\cong M_{2^n}$, for all
$n\ge 1$, and whose Bratteli diagram looks like
\[
\xymatrix{ {\bullet} \ar@<3pt>[r]\ar@<-3pt>[r] & {\bullet}
\ar@<3pt>[r]\ar@<-3pt>[r] & {\bullet} \ar@<3pt>[r]\ar@<-3pt>[r] &
\cdots }.
\]
From Lemma \ref{dimgrplemma}, we know that the dimension group of
$\mathfrak{A}$ is equal to the direct limit of the scaled ordered
groups
\[
\left(\widetilde{C}(U_n,\mathbb{Z}),
\widetilde{C}(U_n,\mathbb{Z}^+), 2^n\chi_{_{U_n}}\right)
\]
where, in the case of the CAR algebra, $U_n=B(1,n)$ in the
notation of \S \ref{section2}.

For each $n\ge 1$, we define the homomorphism
\[
\phi_n:\widetilde{C}(U_n,\mathbb{Z})\to \mathbb{Z}[1/2],
\]
where $\mathbb{Z}[1/2]$ represents the diadic rationals, by
$\phi_n\left( \alpha\chi_{_{B(1,n)}}\right)=\alpha 2^{-n}$. It can
easily be shown that the diagram
\[
\xymatrix{ \widetilde{C}(U_m,\mathbb{Z}) \ar[rr]^{\phi_{(m,n)*}}
\ar[drr]_{\phi_m} & & \widetilde{C}(U_n,\mathbb{Z}) \ar[d]^{\phi_n}\\
& & \mathbb{Z}[1/2]}
\]
commutes.  It follows that
\[
\lim_{\longrightarrow} C(U_n,\mathbb{Z})\cong
\lim_{\longrightarrow} \widetilde{C}(U_n,\mathbb{Z})\cong
\mathbb{Z}[1/2].
\]

We also see that the positive cone, cf.~\cite{davidson}, is given
as
\[
K_0^+(\mathfrak{A})=\bigcup_{n\ge 1}\phi_n\left(
\widetilde{C}(U_n,\mathbb{Z}^+)\right) =\mathbb{Z}^+[1/2],
\]
with $u=1$, where $u$ is the order unit of $\mathbb{Z}[1/2]$.
Thus, we have the dimension group of the CAR algebra
$\mathfrak{A}$ given by
\[
K_0(\mathfrak{A})=\left(\mathbb{Z}[1/2],\mathbb{Z}^+[1/2],1\right).
\]
We note that this example also appears in~\cite{davidson}, and our
calculations here could also be obtained using the results of
Poon~\cite[Corollary 3.3]{poon}.

\end{example}

\begin{example}
\label{continuousoncantor}

The AF algebra with Bratelli diagram
\[
\xymatrix{ & & & {\bullet} \ar[dll]
 \ar[drr] & & &\\
  & {\bullet} \ar[dl] \ar[dr]
   & & & & {\bullet}
  \ar[dl] \ar[dr] &\\
   & & & {\vdots} & & & }
\]
is simply the set of all continuous complex-valued functions on
the Cantor set $X$.  This algebra is mentioned at the beginning of
\S \ref{section3}, and from the discussion given there one has
that
\[
K_0(C(X))=\left(
C(X,\mathbb{Z}),C(X,\mathbb{Z}^+),\chi_{_{X}}\right).
\]

\end{example}

\begin{example}

We now consider an AF algebra which is a hybrid of the two
previous examples. In particular, let $\mathfrak{A}$ be the AF
algebra with Bratteli diagram
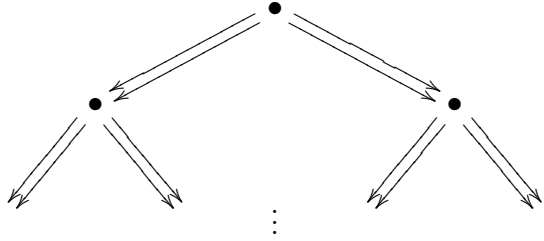
\[
\xymatrix{ & & & {\bullet} \ar@<2pt>[dll]\ar@<-2pt>[dll]
 \ar@<2pt>[drr]\ar@<-2pt>[drr] & & &\\
  & {\bullet} \ar@<2pt>[dl]\ar@<-2pt>[dl] \ar@<2pt>[dr]\ar@<-2pt>[dr]
   & & & & {\bullet}
  \ar@<2pt>[dl]\ar@<-2pt>[dl] \ar@<2pt>[dr]\ar@<-2pt>[dr] &\\
   & & & {\vdots} & & & }
\]
By Lemma \ref{dimgrplemma} we know that $\displaystyle
K_0(\mathfrak{A})=\lim_{\longrightarrow}
\widetilde{C}(U_n,\mathbb{Z})$ where
$U_n=\bigcup_{i=1}^{2^n}B(i,n)$, for all $n\ge 1$, can be thought
of diagrammatically as the union of the sets at level $n$ in the
Bratelli diagram. For $n\ge 1$ given and $\sum_{i=1}^{2^n}
\alpha_i\chi_{_{B(i,n)}}\in\widetilde{C}(U_n,\mathbb{Z})$, we see
that
\begin{eqnarray*}
\phi_{(n,n+1)\ast}\left(\sum_{i=1}^{2^n}\alpha_i\chi_{_{B(i,n)}}\right)
&=&\sum_{i=1}^{2^n}\alpha_i\phi_{(n,n+1)\ast}\left(\chi_{_{B(i,n)}}\right)\\
&=&\sum_{i=1}^{2^n}\alpha_i\left( 2\chi_{_{B(i_1,n+1)}}+2\chi_{_{
B(i_2,n+1)}}\right)\\
&=&2\sum_{i=1}^{2^n}\alpha_i\chi_{_{\left[B(i_1,n+1)\cup
B(i_2,n+1)\right]}}
\end{eqnarray*}
where $B(i_1,n+1)$ and $B(i_2,n+1)$ are exactly those sets in the
collection $\bigvee_{r=1}^{m_{n+1}}B(r,n+1)$ which lie in
$B(i,n)$.

In general, for $m\le n$ we have
\[
\phi_{(m,n)\ast}\left( \sum_{i=1}^{2^m}\alpha_i
\chi_{_{B(i,m)}}\right) =2^{n-m}
\sum_{i=1}^{2^m}\alpha_i\chi_{_{\left[\bigcup_{j=1}^{t(i)}B(r_j,n)\right]}}
\]
where $B(r_1,n),\ldots,B(r_{t(i)},n)$ are exactly those sets in
$\bigvee_{r=1}^{m_n}B(r,n)$ which lie in $B(i,m)$.

Now, let $X_{min}=\bigcap_{n=1}^\infty U_n$, which can be shown to
be the Cantor set. Indeed, the diagram for $X_{min}$ is that of
Example \ref{continuousoncantor}. We then define a homomorphism
\[
\phi_n:\widetilde{C} (U_n,\mathbb{Z})\to
C(X_{min},\mathbb{Z}[1/2]),
\]
for all $n\ge 1$ and where we give $\mathbb{Z}[1/2]$ the discrete
topology, by
\[
\phi_n\left( \sum_{i=1}^{2^n}\alpha_i\chi_{_{B(i,n)}}\right)
=2^{-n}\sum_{i=1}^{2^n}\alpha_i\chi_{_{\left[B(i,n)\cap
X_{min}\right]}}.
\]
We see that for any $m\le n$,
\begin{eqnarray*}
\phi_n\circ\phi_{(m,n)\ast}
\left(\sum_{i=1}^{2^m}\alpha_i\chi_{_{B(i,m)}}\right)&=&
\phi_n\left(2^{n-m}\sum_{i=1}^{2^m}\alpha_i
\chi_{_{\left[\bigcup_{j=1}^{t(i)}B(r_j,n)\right]}}\right)\\
&=&2^{-m}\sum_{i=1}^{2^m}\alpha_i\chi_{_{\left[\bigcup_{j=1}^{t(i)}
B(r_j,n)\right]\cap X_{min}}}\\
&=&2^{-m}\sum_{i=1}^{2^m}\alpha_i \chi_{_{\left[B(i,m)\cap X_{min}\right]}}\\
&=&\phi_m\left( \sum_{i=1}^{2^m}\alpha_i\chi_{_{B(i,m)}}\right).
\end{eqnarray*}
Hence, the diagram
\[
\xymatrix{ \widetilde{C}(U_m,\mathbb{Z}) \ar[r]^{\phi_{(m,n)*}}
\ar[dr]_{\phi_m} & \widetilde{C}(U_n,\mathbb{Z}) \ar[d]^{\phi_n}\\
& C(X_{min},\mathbb{Z}[1/2])}
\]
commutes.  Therefore, $\displaystyle \lim_{\longrightarrow}
\widetilde{C}(U_n,\mathbb{Z})\subset C(X_{min},\mathbb{Z}[1/2])$.
However, if we let $f\in C(X_{min},\mathbb{Z}[1/2])$, then, by the
compactness of $X_{min}$, $f$ is a simple function, and therefore
an element of $\displaystyle \lim_{\longrightarrow}\widetilde{C}
(U_n,\mathbb{Z})$.  Hence, by Lemma \ref{dimgrplemma},
\[
\lim_{\longrightarrow}C(U_n,\mathbb{Z})\cong C(X_{min},
\mathbb{Z}[1/2]).
\]
In addition to this,
\[
K_0^+(\mathfrak{A})=\bigcup_{n\ge
1}\phi_n\left(\widetilde{C}(U_n,\mathbb{Z}^+)\right)=\{f\in
C(X_{min},\mathbb{Z}[1/2]): f\ge 0\}
\]
and the order unit $u$ is given by $u=\chi_{_{X_{min}}}$.  Thus,
\[
K_0(\mathfrak{A})=\left(
C(X_{min},\mathbb{Z}[1/2]),C(X_{min},\mathbb{Z}[1/2]^+),\chi_{_{X_{min}}}\right).
\]

\end{example}

\begin{example}
\label{gicar}

Consider the GICAR algebra $\mathfrak{A}$, whose Bratteli diagram
is given by
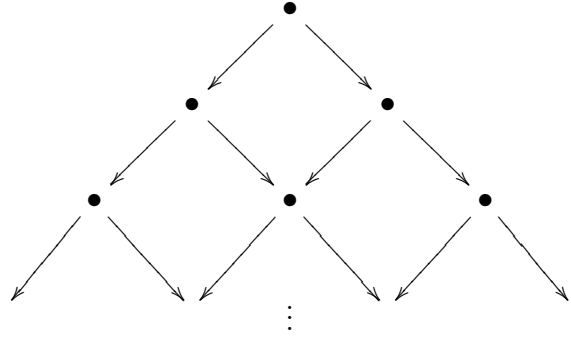
\[
\xymatrix{ & & & {\bullet} \ar[dl] \ar[dr] & & &\\
  & & {\bullet} \ar[dl] \ar[dr] & & {\bullet} \ar[dl] \ar[dr] & &\\
  & {\bullet} \ar[dl] \ar[dr] & & {\bullet} \ar[dl] \ar[dr] & &
  {\bullet} \ar[dl] \ar[dr] &\\
  & & & {\vdots} & & & }
\]
We would like to use Lemma \ref{dimgrplemma} to calculate the
$K_0$ group of $\mathfrak{A}$.

First, from \S \ref{section2}, we know that we can use, as our
decreasing sequence $\{U_n\}_{n=0}^\infty$ of sets, the collection
where $U_n=\cup_{r=1}^{n+1} B(r,n)$, for all $n\ge 1$, with the
clopen sets $\{B(r,n)\}_{r=1}^{n+1}$ corresponding to the upper
left matrix units in the subalgebra $\mathfrak{A}_n$.

The direct system
\[
\widetilde{C}(U_1,\mathbb{Z})\stackrel{\phi_{(1,2)\ast}}{\longrightarrow}
\widetilde{C}(U_2,\mathbb{Z})\stackrel{\phi_{(2,3)\ast}}{\longrightarrow}
\cdots
\]
that results from Theorem \ref{dimgroup} and Lemma
\ref{dimgrplemma} is then such that the action of
$\phi_{(n,n+1)\ast}$ can be realized by the matrix
\[
A=\left[
\begin{array}{ccccc}
1 &   & & &\\
1 & 1 & & &\\
  & 1 & & &\\
  & & & \ddots &\\
  &  & & &1\\
  &  & & &1\\
\end{array}
\right]\in M_{n+2,n+1},
\]
where all unspecified entries are zero and the function
$\sum_{r=1}^{n+1}\alpha_r\chi_{_{B(r,n)}}\in\widetilde{C}(
U_n,\mathbb{Z})$ is viewed as the vector $(\alpha_1,\ldots,
\alpha_{n+1})^T\in \mathbb{Z}^{n+1}$.  We now define, by adding a
column to this matrix, the new matrix
\[
A_{n,n+1}=\left[
\begin{array}{cccccc}
1 & & & & &\\
1 &1& & & &\\
 & 1 & & & &\\
 & & \ddots & & &\\
  & & & 1 &1 &\\
  & & & &1 & 1\\
\end{array}
\right]\in M_{n+2},
\]
where, as above, all unspecified entries are zero. This matrix is
lower triangular, and therefore invertible.  So, for all $n\ge 1$,
let
\[
A_n=\left[A_{0,1}^{-1}\oplus I_{n-1}\right]\left[
A_{1,2}^{-1}\oplus I_{n-2}\right] \cdots \left[ A_{n-2,n-1}^{-1}
\oplus I_1\right]A_{n-1,n}^{-1}\in M_{n+1},
\]
and $R_n:\mathbb{Z}^{n+1}\to C(X_{min},\mathbb{Z})$ be given by
\[
R_n(\alpha_1,\ldots,\alpha_{n+1})^T=\sum_{r=1}^{n+1} \alpha_r
\chi_{_{B(r,r-1)\cap X_{min}}},
\]
where, as before, $X_{min}=\bigcap_{n=0}^\infty U_n$.

\begin{remark}

The set $X_{min}$ can be given a concrete realization as the set
of all infinite paths in the diagram
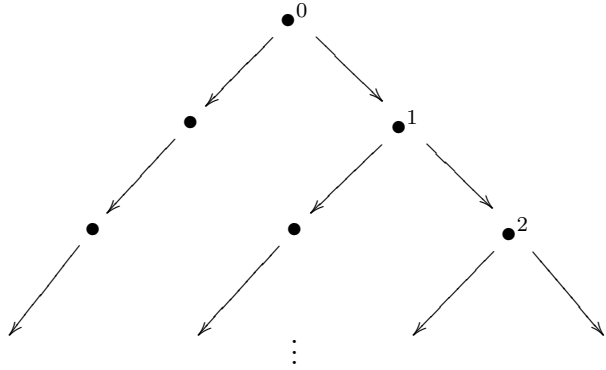
\[
\xymatrix{ & & & {\bullet}^0 \ar[dl] \ar[dr] & & &\\
  & & {\bullet} \ar[dl] & & {\bullet}^1 \ar[dl] \ar[dr] & &\\
  & {\bullet} \ar[dl]  & & {\bullet} \ar[dl]  & &
  {\bullet}^2 \ar[dl] \ar[dr] &\\
  & & & {\vdots} & & & }
\]
which is homeomorphic to the set $\{0,1,1/2,1/3,1/4,\ldots\}$ with
the relative topology it inherits as a subset of $\mathbb{R}$.
Then, the set $B(r,r-1)$, for any $r\ge 1$, is the set of all
paths emanating from the vertex with label $r-1$.

\end{remark}

Finally, let $\phi_n:\widetilde{C}(U_n,\mathbb{Z})\to
C(X_{min},\mathbb{Z})$ be given by $\phi_n=R_n\circ A_n$, for all
$n\ge 1$, where we know $\phi_n$ maps to $C(X_{min},\mathbb{Z})$
since $A_{r,r+1}^{-1}$ has integer entries for all $r\ge 0$.  We
would like to show that the diagram
\[
\xymatrix{ \widetilde{C}(U_n,\mathbb{Z}) \ar[rr]^{\phi_{(n,n+1)*}}
\ar[drr]_{\phi_n} & & \widetilde{C}(U_{n+1},\mathbb{Z}) \ar[d]^{\phi_{n+1}}\\
& & C(X_{min},\mathbb{Z})}
\]
commutes.

Let $(\alpha_1,\ldots,\alpha_{n+1})^T\in\widetilde{C}
(U_n,\mathbb{Z})$.  Then,
\begin{eqnarray*}
\lefteqn{\phi_{n+1}\circ\phi_{(n,n+1)\ast}(\alpha_1,\ldots,
\alpha_{n+1})^T}\\
&=& \phi_{n+1}\circ A_{n,n+1}(\alpha_1,\ldots,\alpha_{n+1},0)^T\\
&=& R_{n+1}\circ A_{n+1}\circ A_{n,n+1}
(\alpha_1,\ldots,\alpha_{n+1},0)^T \\
&=& R_{n+1}\circ [A_{0,1}^{-1}\oplus I_n]\cdots [A_{n-1,n}^{-1}
\oplus I_1](\alpha_1,\ldots,\alpha_{n+1},0)^T.
\end{eqnarray*}
If we define $(\beta_1^{n-1},\ldots,\beta_{n+1}^{n-1})^T=
A_{n-1,n}^{-1}(\alpha_1,\ldots,\alpha_{n+1})^T$ and, in general,
for $2\le i<n$,
\[
(\beta_1^{n-i},\ldots,\beta_{n-i+2}^{n-i})^T= A_{n-i,n-i+1}^{-1}
(\beta_1^{n-i+1},\ldots,\beta_{n-i+2}^{n-i+1})^T,
\]
we see that
\begin{eqnarray*}
\lefteqn{\phi_{n+1}\circ\phi_{(n,n+1)\ast}(\alpha_1,
\ldots,\alpha_{n+1})^T}\\
&=& R_{n+1}\circ [A_{0,1}^{-1}\oplus I_n]\cdots [A_{n-2,n-1}^{-1}
\oplus I_2](\beta_1^{n-1},\ldots,\beta_{n+1}^{n-1},0)^T\\
&=& R_{n+1}\circ [A_{0,1}^{-1}\oplus I_n]\cdots [A_{n-3,n-2}^{-1}
\oplus I_3](\beta_1^{n-2},\ldots,\beta_n^{n-2},\beta_{n+1}^{n-1},
0)^T\\
&=& \cdots =R_{n+1}\circ [A_{0,1}^{-1}\oplus I_n](
\beta_1^1,\beta_2^1,\beta_3^1,\ldots,\beta_n^{n-2},\beta_{n+1}^{n-1},
0)^T\\
&=& R_{n+1}(\beta_1^0,\beta_2^0,\beta_3^1,\ldots, \beta_n^{n-2},
\beta_{n+1}^{n-1},0)^T\\
&=& \beta_1^0\chi_{_{B(1,0)\cap X_{min}}}+ \sum_{r=2}^{n+1}
\beta_r^{r-2}\chi_{_{B(r,r-1)\cap X_{min}}}.
\end{eqnarray*}
However, we see from the following calculation that
\begin{eqnarray*}
\lefteqn{\phi_n(\alpha_1,\ldots,\alpha_{n+1})^T}\\
&=& R_n\circ A_n
(\alpha_1,\ldots,\alpha_{n+1})^T\\
&=& R_n\circ [A_{0,1}^{-1}\oplus I_{n-1}]\cdots [A_{n-2,n-1}^{-1}
\oplus I_1]A_{n-1,n}^{-1}(\alpha_1,\ldots,\alpha_{n+1})^T\\
&=& R_n\circ [A_{0,1}^{-1}\oplus I_{n-1}]\cdots [A_{n-2,n-1}^{-1}
\oplus I_1](\beta_1^{n-1},\ldots,\beta_{n+1}^{n-1})^T\\
&=&\cdots =R_n(\beta_1^0,\beta_2^0,\beta_3^1,\ldots,
\beta_n^{n-2},\beta_{n+1}^{n-1})^T\\
&=& \beta_1^0\chi_{_{B(1,0)\cap X_{min}}}+ \sum_{r=2}^{n+1}
\beta_r^{r-2}\chi_{_{B(r,r-1)\cap X_{min}}}.
\end{eqnarray*}
Thus, $\phi_n=\phi_{n+1}\circ\phi_{(n,n+1)\ast}$, and we conclude
that the diagram commutes.

Now, in constructing the algebra $\mathfrak{A}$, by choosing the
appropriate embeddings between the finite-dimensional subalgebras,
we may assume, without a loss of generality, that the sequence
$\{U_n\}_{n=0}^\infty$ is such that $B(r,n)\subset B(r,n-1)$, for
all $1\le r\le n$, and that $B(n+1,n)\subset B(n,n-1)$.  For $f\in
C(X_{min},\mathbb{Z})$, by the compactness of $X_{min}$, there
exists an $n\ge 1$ and integers $\alpha_1,\ldots,\alpha_{n+1}$
such that
\[
f=\sum_{r=1}^{n+1} \alpha_r\chi_{_{B(r,n)\cap X_{min}}}.
\]
We would like to show that $f$ is an element of the image of
$\phi_n$.

Let $1\le r\le n$.  Then, by assumption,
\[
B(r,n)\subset B(r,n-1)\subset\cdots\subset B(r,r-1).
\]
So, in particular, if $r<l\le n+1$, then $B(r,n)\subset B(r,l-1)$.
Because the collection $\{B(i,l-1)\}_{i=1}^l$ is pairwise
disjoint, it follows that $B(r,n)\cap B(l,l-1)=\emptyset$.  We
also observe that our assumptions imply
\[
B(r,n)\subset B(r,r-1)\subset B(r-1,r-2)\subset\cdots\subset
B(1,0).
\]

Now, define the function
\[
f_{n+1}=\alpha_1\chi_{_{B(1,0)\cap X_{min}}}+ \sum_{i=1}^n
(\alpha_{i+1}-\alpha_i)\chi_{_{B(i+1,i)\cap X_{min}}}.
\]
Take $x\in X_{min}$, and suppose $x\in B(r,n)$, $1\le r\le n$.
Then,
\[
f_{n+1}(x)=\alpha_1\chi_{_{B(1,0)\cap X_{min}}}(x)
+\sum_{i=1}^{r-1} (\alpha_{i+1}-\alpha_i)\chi_{_{B(i+1,i)\cap
X_{min}}}(x).
\]
But, $B(r,n)\subset B(r,r-1)\subset\cdots\subset B(1,0)$, from
which it follows that
\[
f_{n+1}(x)=\alpha_1+\sum_{i=1}^{r-1}(\alpha_{i+1}-\alpha_i)=
\alpha_r.
\]
Thus, $f(x)=\alpha_r=f_{n+1}(x)$, whenever $x\in B(r,n)$, $1\le
r\le n$.  If $x\in B(n+1,n)$ then $f(x)=\alpha_{n+1}$ and
\begin{eqnarray*}
f_{n+1}(x)&=&\alpha_1\chi_{_{B(1,0)\cap X_{min}}}(x)+ \sum_{i=1}^n
(\alpha_{i+1}-\alpha_i)\chi_{_{B(i+1,i)\cap
X_{min}}}(x)\\
&=&\alpha_{ n+1}
\end{eqnarray*}
since $B(n+1,n)\subset B(n,n-1)\subset\cdots\subset B(1,0)$.
Therefore, $f=f_{n+1}$ as elements of $C(X_{min},\mathbb{Z})$. So,
for all $f\in C(X_{min},\mathbb{Z})$ there exist integers
$\beta_1,\ldots,\beta_{n+1}$ such that
\[
f=\sum_{r=1}^{n+1}\beta_r\chi_{_{B(r,r-1)\cap X_{min}}}.
\]

\begin{remark}

The continuous functions $C(X_{min},\mathbb{Z})$ can be identified
with the set of all infinite sequences which are eventually zero.

\end{remark}

Now, since $A_n$ is invertible and $A_n^{-1}$ has integer entries,
\[
A_n^{-1}(\beta_1,\ldots,\beta_{n+1})^T\in\mathbb{Z}^{ n+1},
\]
and $\phi_n(A_n^{-1}(\beta_1,\ldots,\beta_{n+1})^T)=f$.  Thus, we
have
\[
\bigcup_{n\ge 1}\phi_n(\widetilde{C} (U_n,\mathbb{Z}))=
C(X_{min},\mathbb{Z}).
\]
The injectivity of the maps $\phi_n$ then implies that
\[
\lim_{\longrightarrow}\widetilde{C}(U_n,\mathbb{Z}) =\bigcup_{n\ge
1}\phi_n(\widetilde{C}(U_n,\mathbb{Z}))=C(X_{min},\mathbb{Z}).
\]
It also follows from Theorem \ref{dimgroup} that the order unit
$u$ is the function $\chi_{_{X_{min}}}$. To obtain a complete
description of the dimension group of $\mathfrak{A}$, it remains
to describe the positive cone $K_0^+(\mathfrak{A})$.

To begin, as a notational convenience we will let the collection
$\{e_r(n+1)\}_{r=1}^{n+1}$ represent the standard basis vectors
for $\mathbb{Z}^{n+1}$ for each $n\ge 1$.  Our first step will be
to describe the way the matrix $A_n$ acts on this collection. In
what follows we will use $\displaystyle \zerr{a}{b}$ to denote the
binomial coefficient $\displaystyle \frac{a!}{b!(a-b)!}$.

\begin{lemma}

For each $n\ge 1$ and $1\le r\le n+1$, we have
\[
A_n(e_r(n+1))=\sum_{j=r}^{n+1}(-1)^{j-r}\zerr{n+1-r}{j-r}
e_j(n+1).
\]

\end{lemma}

\begin{proof}

We proceed by induction on $n$.  For the induction basis, it is
easy to verify that the formula holds when $n=1$ and $1\le r\le
2$.  So, we suppose the formula is valid for $n\ge 1$ and all
$1\le r\le n+1$.

Let $1\le i\le n+2$ be arbitrary.  We see that
\begin{eqnarray*}
\lefteqn{A_{n+1}(e_i(n+2))}\\
&=&\left[\begin{array}{cc} A_n & 0\\ 0 &
1\end{array}\right] A_{n,n+1}^{-1}e_i(n+2)\\
&=&\left[\begin{array}{cc} A_n & 0\\ 0& 1\end{array}\right]
\left(\sum_{r=i}^{n+2}(-1)^{r-i}e_r(n+2)\right)\\
&=&\sum_{r=i}^{n+2}(-1)^{r-i}\left[ \begin{array}{cc} A_n & 0\\ 0
& 1\end{array}\right] e_r(n+2)\\
&=&\sum_{r=i}^{n+1}(-1)^{r-i}\left[ \begin{array}{cc} A_n & 0\\ 0
& 1\end{array}\right] e_r(n+2)+(-1)^{n-i}e_{n+2}(n+2)\\
&=&\sum_{r=i}^{n+1}(-1)^{r-i}\zerr{A_ne_r(n+1)}{0}+(-1)^{n-i}
e_{n+2}(n+2).
\end{eqnarray*}
By the induction hypothesis we then have
\begin{eqnarray*}
\lefteqn{A_{n+1}(e_i(n+2))}\\
&=&\sum_{r=i}^{n+1}\sum_{j=r}^{n+1} (-1)^{j-i}\zerr{ n+1-r}{j-r}
e_j(n+2)+(-1)^{n-i}e_{n+2}(n+2).
\end{eqnarray*}
If, in this double sum, we switch the order of summation, this
becomes
\begin{eqnarray*}
\lefteqn{A_{n+1}(e_i(n+2))}\\
&=&\sum_{j=i}^{n+1}\sum_{r=i}^j(-1)^{j-i}\zerr{n+1-r}{j-r}
e_j(n+2)+(-1)^{n-i}e_{n+2}(n+2)\\
&=&\sum_{j=i}^{n+1}(-1)^{j-i}\left[\sum_{r=i}^j \zerr{n+1-r}{j-r}
\right] e_j(n+2)+(-1)^{n-i}e_{n+2}(n+2).
\end{eqnarray*}

Here, we make note of the binomial coefficient identity
\[
\sum_{k=0}^n \zerr{m+k}{k}=\zerr{m+n+1}{n},
\]
which is valid for any nonnegative integers $n$ and $m$. By then
letting $k=j-r$ and $m=n+1-j$, we see that we can write
\begin{eqnarray*}
\sum_{r=i}^j\zerr{n+1-r}{j-r}&=& \sum_{k=j-i}^0\zerr{n+1+k-j}{k}\\
&=&\sum_{k=0}^{j-i}\zerr{m+k}{k}\\
&=&\zerr{m+j-i+1}{j-i}\\
&=&\zerr{n-i+2}{j-i}.
\end{eqnarray*}
Therefore,
\begin{eqnarray*}
\lefteqn{A_{n+1}(e_i(n+2))}\\
&=&\sum_{j=i}^{n+1}(-1)^{j-i}\zerr{n-i+2}{j-i}e_j(n+2) +(-1)^{n-i}
e_{n+2}(n+2)\\
&=&\sum_{j=i}^{n+2}(-1)^{j-i}\zerr{n-i+2}{j-i}e_j(n+2).
\end{eqnarray*}
Hence, by induction, the formula holds for all $n\ge 1$ and $1\le
r\le n+1$, completing the proof.

\end{proof}

This gives us the tools we need to describe $K_0^+(\mathfrak{A})$.
However, we first establish some notation.  Let $n\ge 1$ be fixed
and define
\[
P_n=\left\{\sum_{r=1}^{n+1}\beta_r\chi_{_{B(r,r-1)\cap X_{min}}}
:\beta_r\in\mathbb{Z} \mbox{ for all } 1\le r\le n+1\right\}.
\]
We then let $P_n^+$ be the subset of $P_n$ for which
\[
\sum_{l=0}^k \zerr{n-l}{k-l}\beta_{l+1}\ge 0,
\]
for all $k$ such that $0\le k\le n$.  We then have the following
theorem describing $K_0^+$.

\begin{theorem}

For the GICAR algebra $\mathfrak{A}$, $\displaystyle
K_0^+(\mathfrak{A})=\bigcup_{n\ge 1}P_n^+$.

\end{theorem}

\begin{proof}

To begin, let $f\in K_0^+(\mathfrak{A})=\bigcup_{n\ge 1}
\phi_n\left( \widetilde{C} (U_n,\mathbb{Z}^+)\right)$.  So, there
exists an $n\ge 1$ and
$(\alpha_1,\ldots,\alpha_{n+1})^T\in\mathbb{Z}_+^{n+1}$ such that
$f=\phi_n(g)$ where
\[
g=\sum_{r=1}^{n+1}\alpha_r\chi_{_{B(r,n)}}.
\]
We will analyze the action of $\phi_n$ on $g$ by considering the
action of the matrix $A_n$ on the vector $(\alpha_1,\ldots,
\alpha_{n+1})^T$.  So, we apply our previous lemma to get
\begin{eqnarray*}
A_n(\alpha_1,\ldots,\alpha_{n+1})^T&=& A_n\left( \sum_{r=1}^{n+1}
\alpha_r e_r(n+1)\right)\\
&=&\sum_{r=1}^{n+1}\alpha_r A_n(e_r(n+1))\\
&=&\sum_{r=1}^{n+1}\alpha_r \sum_{j=r}^{n+1}(-1)^{j-r}
\zerr{n+1-r}{j-r} e_j(n+1)\\
&=&\sum_{r=1}^{n+1}\sum_{j=r}^{n+1}(-1)^{j-r}\alpha_r
\zerr{n+1-r}{j-r} e_j(n+1).
\end{eqnarray*}
By switching the order of summation this becomes
\[
A_n(\alpha_1,\ldots,\alpha_{n+1})^T= \sum_{j=1}^{n+1} \sum_{r=1}^j
(-1)^{j-r}\alpha_r\zerr{n+1-r}{j-r}e_j(n+1)
\]
and so we have
\begin{eqnarray*}
\phi_n(g)&=& R_n\circ A_n(\alpha_1,\ldots,\alpha_{n+1})^T\\
&=& R_n\left(\sum_{j=1}^{n+1}\sum_{r=1}^j (-1)^{j-r}\alpha_r
\zerr{n+1-r}{j-r}e_j(n+1)\right)\\
&=& \sum_{j=1}^{n+1}R_n\left( \sum_{r=1}^j (-1)^{j-r}\alpha_r
\zerr{n+1-r}{j-r}e_j(n+1)\right)\\
&=&\sum_{j=1}^{n+1}\sum_{r=1}^j(-1)^{j-r}\alpha_r
\zerr{n+1-r}{j-r}\chi_{_{B(j,j-1)\cap X_{min}}}.
\end{eqnarray*}
Thus, we need to show that, for all $0\le k\le n$,
\[
S_k=\sum_{l=0}^k\zerr{n-l}{k-l}\sum_{r=1}^{l+1} (-1)^{l+1-r}
\alpha_r\zerr{n+1-r}{l+1-r}\ge 0.
\]

If we now reverse the order of summation, $S_k$ can be written as
\[
S_k=\sum_{r=1}^{k+1}\sum_{l=r-1}^k (-1)^{l+1-r}\alpha_r
\zerr{n-l}{k-l}\zerr{n+1-r}{l+1-r}.
\]
By making use of the identity
\[
\zerr{n-l}{k-l}\zerr{n+1-r}{l+1-r}=\zerr{n+1-r}{k+1-r}
\zerr{k+1-r}{l+1-r},
\]
we are able to write
\begin{eqnarray*}
S_k&=& \sum_{r=1}^{k+1}\sum_{l=r-1}^k(-1)^{l+1-r}\alpha_r
\zerr{n+1-r}{k+1-r}\zerr{k+1-r}{l+1-r}\\
&=&\sum_{r=1}^{k+1}\alpha_r\zerr{n+1-r}{k+1-r} \sum_{l=r-1}^k
(-1)^{l+1-r}\zerr{k+1-r}{l+1-r}.
\end{eqnarray*}
By a well-known binomial coefficient identity,
\[
\sum_{l=r-1}^k(-1)^{l+1-r}\zerr{k+1-r}{l+1-r} =0,
\]
provided $r-1<k$.  Thus, $S_k=\alpha_{k+1}\ge 0$.  Therefore,
$f\in P_n^+$, and we have established containment in one
direction.

To prove the opposite containment, let $n\ge 1$ be given and
suppose the function $\sum_{r=1}^{n+1} \beta_r\chi_{_{
B(r,r-1)\cap X_{min}}}$ is an element of $P_n^+$.  It follows that
for all $0\le k\le n$ we have $\alpha_k\ge 0$ where
\[
\alpha_k=\sum_{l=0}^k\zerr{n-l}{k-l}\beta_{l+1}.
\]
Hence, the function $\sum_{r=1}^{n+1}\alpha_{r-1}\chi_{_{
B(r,n)}}$ is an element of $\widetilde{C}(U_n,\mathbb{Z}^+)$, and
therefore
\[
\phi_n\left(\sum_{r=1}^{n+1}\alpha_{r-1}\chi_{_{B(r,n)}}\right)
\in K_0^+(\mathfrak{A}).
\]
So, we will have proven the desired result if we can show that
\[
\phi_n\left(\sum_{r=1}^{n+1}\alpha_{r-1}\chi_{_{B(r,n)}}\right)
=\sum_{r=1}^{n+1}\beta_r\chi_{_{B(r,r-1)\cap X_{min}}}.
\]

In order to do this, we first notice that
\begin{eqnarray*}
A_n(\alpha_0,\ldots,\alpha_n)^T&=& A_n\left(\sum_{r=1}^{n+1}
\alpha_{r-1} e_r(n+1)\right)\\
&=&A_n\left(\sum_{r=1}^{n+1}\sum_{l=0}^{r-1} \zerr{n-l}{r-1-l}
\beta_{l+1}e_r(n+1)\right)\\
&=&\sum_{r=1}^{n+1}\sum_{l=0}^{r-1}\zerr{n-l}{r-1-l}\beta_{l+1}
A_n(e_r(n+1)),
\end{eqnarray*}
which, if we apply our previous lemma becomes
\[
\sum_{r=1}^{n+1}\sum_{l=0}^{r-1}\sum_{j=r}^{n+1}(-1)^{j-r}
\beta_{l+1}\zerr{n-l}{r-1-l}\zerr{n+1-r}{j-r}e_j(n+1).
\]
By switching the order of summation, we can in turn write this sum
as
\[
\sum_{j=1}^{n+1} \sum_{r=1}^j
\sum_{l=0}^{r-1}(-1)^{j-r}\beta_{l+1} \zerr{n-l}{r-1-l}
\zerr{n+1-r}{j-r}e_j(n+1).
\]
So, our task has become one of showing that, for each $1\le j\le
n+1$, if
\[
S_j=\sum_{r=1}^j\sum_{l=0}^{r-1}(-1)^{j-r}\beta_{l+1}
\zerr{n-l}{r-1-l}\zerr{n+1-r}{j-r}
\]
then $S_j=\beta_j$.

To establish this, we first switch the order of summation in $S_j$
to obtain
\[
S_j=\sum_{l=0}^{j-1}\sum_{r=l+1}^j(-1)^{j-r}\beta_{l+1}
\zerr{n-l}{r-1-l}\zerr{n+1-r}{j-r},
\]
and notice that if $0\le l<j-1$ then
\begin{eqnarray*}
\lefteqn{\sum_{r=l+1}^j(-1)^{j-r}\zerr{n-l}{r-1-l}\zerr{n+1-r}{j-r}}\\
&=&\sum_{r=l+1}^j(-1)^{j-r}\zerr{j-l-1}{j-r} \frac{(n-l)!}{
(n+1-j)!(j-l-1)!}\\
&=& \frac{(n-l)!}{(n+1-j)!(j-l-1)!} \sum_{r=l+1}^j(-1)^{j-r}
\zerr{j-l-1}{j-r}\\
&=&0
\end{eqnarray*}
since $\sum_{r=l+1}^j(-1)^{j-r}\zerr{j-l-1}{j-r}=0$.  Therefore,
we conclude that $S_j=\beta_j$.  So,
\begin{eqnarray*}
\phi_n(\alpha_0,\ldots,\alpha_n)^T&=& R_n\circ A_n
(\alpha_0,\ldots,\alpha_n)^T\\
&=&R_n\left(\sum_{j=1}^{n+1}\beta_je_j(n+1)\right)\\
&=&\sum_{r=1}^{n+1}\beta_r\chi_{_{B(r,r-1)\cap X_{min}}},
\end{eqnarray*}
which completes the proof.

\end{proof}

For comparison, we mention that Davidson's~\cite{davidson}
description of $K_0(\mathfrak{A})$ for the GICAR algebra is
\[
\left(\mathbb{Z}[y], \{f\in\mathbb{Z}[y]:f>0\mbox{ on } (0,1)\},
1\right).
\]

\end{example}

\section{Groupoids from Dimension Groups}
\label{section5}

The examples considered in the previous section lead naturally to
certain additional questions.  One, a partial answer to which we
give below, is the question of how to reverse the process used in
calculating these examples.  The other, for which no answer is
given here, has to do with the common characteristic shared by
each of the examples. Specifically, these $K_0$ groups, as a
result of the methods used to obtain them, all have in common the
fact that they are groups of continuous functions. It is natural
then to wonder if a general result can be obtained which gives the
$K_0$ group of any AF algebra (or possibly only those of a
specific type) as a group of continuous functions.

Returning to the first of the above questions, we note that the
examples of the previous section illustrate a use of the result in
\S \ref{section3} which allows a calculation of
$K_0(C^*(\mathcal{R}))$ for an AF groupoid $\mathcal{R}$ by making
use of the inverse semigroup of partial homeomorphisms which
generate the groupoid.  So, in this sense, we are able to move
directly from an AF groupoid to a dimension group without first
directly constructing the algebra $C^*(\mathcal{R})$.  It is
natural to ask then whether a converse to this operation is
possible. Of course, by using the dimension group to build the
associated AF algebra, cf.~\cite{davidson}, it is possible to
construct the appropriate groupoid.  However, we would like to do
this directly. That is, given a dimension group $G$, how can one
{\it directly} construct the groupoid $\mathcal{R}$ for which
$K_0(C^*(\mathcal{R}))\cong G$?  We give a partial answer to this
question for a certain class of dimension groups in the following
example.

\begin{example}

Consider the dimension group $\displaystyle (G,G^+,u)$ which is
obtained as the direct limit of the sequence $(\mathbb{Z},
\mathbb{Z}^+,u_n)$ of scaled ordered groups.  In this case, it is
not hard to show that
\[
G/\langle u\rangle \cong \lim_{\longrightarrow} \mathbb{Z}/\langle
u_n\rangle,
\]
where $\langle u\rangle$ indicates the subgroup generated by $u$.

Before proceeding, we pause briefly to mention some additional
background.  Given a group $Q$, a character of $Q$ is a
homomorphism of $Q$ into the multiplicative group $\mathbb{T}=\{
z\in\mathbb{C}: |z|=1\}$, the unit circle in the complex plane.
Suppose now that $Q$ is a topological abelian group.  The set of
all continuous characters of $Q$ forms a topological abelian
group, called the character group of $Q$, and will be denoted by
$\widehat{Q}$.  We use~\cite{hewittandross} as a reference on
character groups.

Now, letting $Q=G/\langle u\rangle$ and $Q_n=\mathbb{Z}/\langle
u_n \rangle \cong\mathbb{Z}_{u_n}$, for all $n\ge 1$, we consider
$\widehat{Q}$ and $\widehat{Q}_n$, the character groups of $Q$ and
$Q_n$, respectively.  As an application of the notion of duality,
one can show that the character group of $Q$ is the projective
limit of the character groups of $Q_n$.  In symbols,
\[
\displaystyle \widehat{Q}\cong \lim_{\longleftarrow}\widehat{Q}_n,
\]
and the connecting homomorphisms $\psi_n:\widehat{Q}\to
\widehat{Q}_n$ and $\psi_{n+1,n}:\widehat{Q}_{n+1}\to
\widehat{Q}_n$ are continuous surjections, for all $n\ge 1$, such
that the diagram
\[
\xymatrix{ \widehat{Q}_n & & \widehat{Q}_{n+1} \ar[ll]_{\psi_{n+1,n}}\\
& & \widehat{Q} \ar[u]_{\psi_{n+1}} \ar[ull]^{\psi_n}}
\]
commutes.

Letting $\mbox{id}_1\in\widehat{Q}_1$ be the character which is
identically $1$ on $Q_1$, consider the set $B(1,1)=\psi_1^{-1}
(\mbox{id}_1)$.  By the continuity of $\psi_1$, it is clear that
$B(1,1)$ is a clopen subgroup of $\widehat{Q}$. Hence, for all
characters $\chi_1,\chi_2\in\widehat{Q}$, if the cosets
$\chi_1B(1,1)$ and $\chi_2B(1,1)$ are such that $\chi_1 B(1,1)\cap
\chi_2 B(1,1)\ne\emptyset$ then $\chi_1 B(1,1)=\chi_2 B(1,1)$. The
collection of cosets $\{\chi B(1,1):\chi\in\widehat{Q}\}$ is
clearly an open cover of $\widehat{Q}$, and so by the compactness
of $\widehat{Q}$, there exist $\chi_1^{(1)},\ldots,\chi_n^{(1)}
\in\widehat{Q}$ such that
\[
\widehat{Q}=\bigcup_{i=1}^n\chi_i^{(1)} B(1,1).
\]
This is a disjoint union since we may assume, without a loss of
generality by our above observations, that
$\chi_1^{(1)}=\mbox{id}\in\widehat{Q}$ and $\chi_i^{(1)}
B(1,1)\cap\chi_j^{(1)}B(1,1)=\emptyset$ whenever $i\ne j$.

Now, it is straightforward that $\psi_1(\chi_i^{(1)} B(1,1))$ is a
singleton for each $1\le i\le n$.  If we suppose
$\psi_1(\chi_i^{(1)}B(1,1))$ and $\psi_1(\chi_j^{(1)}B(1,1))$ are
the same element, we find that $i=j$ by our choice of these
characters.  Thus, because $\psi_1$ surjects,
$\psi_1(\widehat{Q})=\widehat{Q}_1=\widehat{\mathbb{Z}}_{ u_1}
=\mathbb{Z}_{u_1}$,  and it must be that $n=u_1$.  So, we obtain
the partition
$$
\bigvee_{i=1}^{u_1}\chi_i^{(1)}B(1,1)
$$
of $\widehat{Q}$.

Now, for each $n\ge 1$, let $B(1,n)=\psi_n^{-1}(\mbox{id}_n)
\subset\widehat{Q}$, which is a clopen subgroup.  We see that if
we take $\chi\in B(1,n+1)$ then
$$
\psi_n(\chi)=\psi_{n+1,n}\circ\psi_{n+1}(\chi) =\psi_{n+1,n}
(\mbox{id}_{n+1})=\mbox{id}_n.
$$
Thus, $\chi\in B(1,n)$, and we see that $\{B(1,n)\}_{n=1}^\infty$
is a nested decreasing sequence of clopen subgroups in
$\widehat{Q}$.

As above, for each $n\ge 1$, there exists characters
$\chi_1^{(n+1)},\ldots,\chi_{s(n+1)}^{(n+1)}\in B(1,n)$ such that
$$
\bigvee_{i=1}^{s(n+1)}\chi_i^{(n+1)}B(1,n+1)
$$
is a partition of $B(1,n)$.

\begin{claim}

$s(n+1)=u_{n+1}/u_n$

\end{claim}

\begin{reason}

Suppose $\bigvee_{i=1}^r\chi_i B(1,n)$ is a partition of
$\widehat{Q}$.  As above, $r=u_n$.  It follows that
$$
\bigvee_{i=1}^{u_n}\bigvee_{j=1}^{s(n+1)}\chi_i\chi_j^{(n+1)}
B(1,n+1)
$$
is a partition of $\widehat{Q}$ and thus that $u_n
s(n+1)=u_{n+1}$. Hence, we have $s(n+1)=u_{n+1}/u_n$.

\end{reason}

For each $n\ge 1$ and with $B(1,0)=\widehat{Q}$, if we now define
the maps $\sigma_{1,s}^{(n)}: B(1,n)\to B(1,n-1)$ by
$\sigma_{1,s}^{(n)}(\chi)=\chi_s^{(n)}\chi$, then clearly the
collection
$$
\bigcup_{n=1}^\infty\bigcup_{s=1}^{s(n)}\{\sigma_{1,s}^{(n)}\}
$$
is a collection of partial homeomorphisms on $\widehat{Q}$ that
satisfy Conditions \ref{conditions} of \S \ref{section2}.  Thus,
for all $n\ge 1$, let $\tau_{(s_1,\ldots,s_n)}^{(n)}: B(1,n)\to
\widehat{Q}$ be given by
$$
\tau_{(s_1,\ldots,s_n)}^{(n)}(\chi)=\left[ \prod_{i=1}^n
\chi_{s_i}^{(i)}\right] \chi,
$$
where $(s_1,\ldots,s_n)\in\prod_{i=1}^n [1,u_i/u_{i-1}]$, with
$u_0=1$, and each of these intervals is a subset of
$\mathbb{Z}^+$. Then, by Corollary \ref{classifycor2}, the inverse
semigroup $S$ generated by the partial homeomorphisms
$\tau_{(s_1,\ldots,s_n)}^{(n)}$ is AF. Furthermore, Lemma
\ref{dimgrplemma} allows us to calculate its dimension group as
the direct limit of the sequence
$$
(\mathbb{Z},u_1)\stackrel{\phi_{(1,2)\ast}}{\longrightarrow}
(\mathbb{Z},u_2)\stackrel{\phi_{(2,3)\ast}}{\longrightarrow}
\cdots
$$
where the connecting maps $\phi_{(n,n+1)\ast}$ are such that
$$
\phi_{(n,n+1)\ast}(x)=\left(u_{n+1}/u_n\right)x.
$$
But then it is clear that with $\mathcal{R}$ the groupoid
generated by the inverse semigroup $S$, we have
\[
(G,G^+,u)\cong (K_0(C^*(\mathcal{R})),K_0(C^*(\mathcal{R})^+),1)
\]
as desired.

\end{example}

\bibliography{authorarxiv}

\end{document}